\documentclass[11pt, twoside]{amsart}

\usepackage[T1]{fontenc}
\usepackage{amssymb,amsmath,amstext}

\newtheorem{theor}{Theorem}[section]
\newtheorem{lem}[theor]{Lemma}
\newtheorem{defin}[theor]{Definition}

\newtheorem{cor}[theor]{Corollary}
\newtheorem{rem}[theor]{Remark}

\newtheorem{fact}[theor]{Fact}

\newtheorem{observation}[theor]{Observation}

\newcommand{\mr}{\mathrm}

\newcommand{\mb}{\mathbf}

\newcommand{\es}{\emptyset}

\newcommand{\dist}{\mathrm{dist}}

\newcommand{\mcC}{\mathcal{C}}

\newcommand{\mcG}{\mathcal{G}}
\newcommand{\mcH}{\mathcal{H}}

\newcommand{\mcM}{\mathcal{M}}

\newcommand{\mcS}{\mathcal{S}}

\newcommand{\mbC}{\mathbf{C}}

\newcommand{\mbG}{\mathbf{G}}
\newcommand{\mbH}{\mathbf{H}}

\newcommand{\mbP}{\mathbf{P}}

\newcommand{\mbX}{\mathbf{X}}

\newcommand{\mbbE}{\mathbb{E}}
\newcommand{\mbbP}{\mathbb{P}}
\newcommand{\mbbN}{\mathbb{N}}

\newcommand{\mbbQ}{\mathbb{Q}}

\title[Random graphs with bounded maximum degree]
{Random graphs with bounded maximum degree: \\ asymptotic structure and a logical limit law}

\author{Vera Koponen}

\address{Vera Koponen, Department of Mathematics, Uppsala University, Box 480,
75106 Uppsala, Sweden.}

\email{vera@math.uu.se}

\begin{document}

\begin{abstract} 
{\Small For any fixed integer $R \geq 2$ we characterise the typical structure
of undirected graphs with vertices $1, \ldots, n$ and maximum degree $R$, as $n$ tends to infinity.
The information is used to prove that such graphs satisfy a labelled limit law for first-order
logic. If $R \geq 5$ then also an unlabelled limit law holds.
}
\end{abstract}

\maketitle

\section{Introduction}\label{introduction}

\noindent
This article is concerned with finite graphs with $n$ vertices labelled $1, \ldots, n$ such that
every vertex has degree at most $R$, where $R$ will be a fixed non-negative integer.
By {\em graph} we always mean finite undirected graph.
The asymptotic structure of such graphs have been studied before, in numerous
articles, in the case when a degree sequence is associated and one considers (only)
graphs with this degree sequence. Two variants of degree sequences have been
considered. In one variant a sequence $\bar{d} = (d_1, \ldots, d_n)$ 
of non-negative integers is given and one studies
the graphs with vertices $1, \ldots, n$ such that, for every $i = 1, \ldots, n$,
the vertex $i$ has degree $d_i$. In the other variant we are given a sequence
$\bar{d} = (d_0, d_1, \ldots, d_R)$ and study the graphs with vertices
$1, \ldots, n$ such that exactly $d_i$ vertices have degree $i$.
Early work on graphs with prescribed degree sequence include 
\cite{BC, Bol80, Wor81, Wor81b, MW84}. 
In particular one can consider degree sequences $(d_1, \ldots, d_n)$
in which all entries are the same, say $r$,
in which case one gets $r$-regular graphs.
This special case has also been studied extensively, for instance as a special case in
the mentioned articles.
See \cite{Wor99} for a survey about random regular graphs.
Suppose that for all sufficiently large integers $n$, 
$\mbG_n$ is a set of graphs with vertices $1, \ldots, n$. We say that
$\mbG_n$ satisfies a labelled logical {\em limit law} if for every
property which can be expressed with a first-order sentence in the language of graphs,
the proportion of graphs in $\mbG_n$ which have this property converges as $n \to \infty$.
If convergence holds in the case when we count graphs only up to non-isomorphism, then
we say that $\mbG_n$ satisfies an unlabelled limit law.
If the number to which the mentioned proportion converges is always 0 or 1,
then one says that a {\em zero-one law} holds.
Lynch \cite{Lyn05} has proved that a labelled logical limit law holds if we consider
an infinite sequence of degree sequences $(\bar{d}_n)_{n=1}^{\infty}$ (of the first kind
mentioned above) which satisfies
certain properties and for each $n$ the graphs with vertices $1, \ldots, n$ that have
the degree sequence $\bar{d}_n = (d_{n,1}, \ldots, d_{n,n})$. 
In particular, the main result in \cite{Lyn05} implies that
$r$-regular graphs satisfy a labelled limit law, but not a zero-one law.
One of the assumed properties of the sequence of degree sequences in \cite{Lyn05} fails in the
instances that turn out to be of interest in this article, so we cannot use \cite{Lyn05}.
Haber and Krivelevich \cite{HK10} have proved zero-one laws, and the absence of a limit law, for
$r(n)$-regular graphs when $r(n)$ is a function
of the form $\alpha n$ or $n^{\alpha}$ where $0 < \alpha < 1$; in the latter case it matters
whether $\alpha$ is rational or not.

The case of random graphs with a fixed bound on the maximum degree, but no other constraint,
has received less attention than the cases discussed so far.
However, papers by Kennedy and Quintas \cite{KQ85} and 
Rucinski and Wormald \cite{RW92, RW97} are exceptions from this neglect and
study the set of graphs with vertices $1, \ldots, n$ such that every vertex has degree at most $R$,
which we denote $\mbG_{n,R}$.
They do this from the point of view of a random process which generates
members of $\mbG_{n,R}$ by randomly adding new edges as long as it is possible without
violating the condition that the maximum degree is $R$.
Kennedy and Quintas \cite{KQ85} stresses the importance of graphs with 
degree limitations from the perspective of models of physical and chemical systems and gives several
references to work concerned with these aspects.
Rucinski and Wormald \cite{RW92} prove that with probability approaching 1 as $n$ tends
to infinity, the random process ends up with a graph with $\lfloor nR/2 \rfloor$ edges.
In \cite{RW97} they prove that if $R=2$, then the number of $k$-cycles ($k \geq 3$) in a graph
produced by the random process has a Poisson distribution.

However, the main motivation of this article for studying $\mbG_{n,R}$ comes from logic and
from understanding the fine structure of $\mcH$-free graphs for certain types of $\mcH$.
From the point of view of model theory in logic, a restriction on the maximum degree
is natural and can be generalised to arbitrary relational structures via the so-called 
Gaifman graph (e.g. \cite{EF}) of the structure. Finite relational structures and
first-order logic correspond closely to relational databases and the standard query 
language SQL, respectively (e.g. \cite{Lib}, Chapter 1.1).
Within the subfield of model theory which studies logical limit laws
one may be interested in better understanding which classes of structures
(e.g. graphs) have a logical limit (or even zero-one) law and
which do not. When a logical limit law exists for a given logical language and collection $\mbC$
of finite structures it follows that, for every property $P$ which can be expressed
in that language, the proportion of structures which satisfy $P$ stabilises as the number
of elements in the structure tends to infinity. 
Such information may be of help when understanding the average performance of algorithms
that take random members of $\mbC$ as input.

Another type of restriction on graphs, studied in graph theory and interesting from a model theoretic
point of view, is to forbid the occurence of subgraphs which are isomorphic to
some specied graph $\mcH$ (or to a member of a specified set of graphs).
The results of this paper are relevant for understanding $\mcH$-free graphs for
certain types of graphs $\mcH$. 
Let us call a vertex $v$ of the graph $\mcH$ {\em colour critical} if the removal
of $v$ from $\mcH$ leaves a graph with lower chromatic number than $\mcH$.
The {\em criticality} of $v$ is the minimal number of edges containing $v$ which
have to be removed in order to reduce the chromatic number.
Hundack, Prömel and Steger \cite{HPS} have proved that if $\mcH$ is a graph with 
chromatic number $l+1$ ($l \geq 2$) and a
colour critical vertex $v$ of criticality $R+1$ and $\mb{Forb}_n(\mcH)$ is the set of $\mcH$-free 
graphs with vertices $1, \ldots, n$, then the proportion of $\mcG \in \mb{Forb}_n(\mcH)$ which
have the following property approaches 1 as $n \to \infty$: the vertex set can be partitioned
into parts $V_1, \ldots, V_l$ such that for every $i$ the induced subgraph with vertex set $V_i$
has maximum degree $R$.
Hence, understanding of the asymptotic fine structure of graphs with maximum degree $R$ helps us 
to understand the asymptotic fine
structure of $\mcH$-free graphs, for certain $\mcH$.
This connection is used in \cite{Kop12} to prove a logical limit law (but not zero-one law) for
$\mcH$-free graphs, when $\mcH$ is a complete $(l+1)$-partite graph ($l \geq 2$) in which one
part is a singleton.

If nothing else is said then we assume that $R \geq 2$ is an integer.
We always consider the uniform probability distribution on $\mbG_{n,R}$.
In Section~\ref{the distribution of degrees} we study the typical degree
distribution of a random $\mcG \in \mbG_{n,R}$ as $n \to \infty$.
We find that, as $n \to \infty$, a random $\mcG \in \mbG_{n,R}$ almost surely has no vertices
with degree less than $R-2$ and, for every $\varepsilon > 0$, 
the number of vertices with degree $R-1$ is almost
surely between $\sqrt{(R - \varepsilon)n}$ and $\sqrt{(R + \varepsilon)n}$,
where {\em almost surely} means that the probability of the event approaches 1 as $n \to \infty$.
The number of vertices with degree $R-2$ has, asymptotically, a Poisson distribution.
These results are summarised in Theorem~\ref{theorem about typical degree distribution}.
To prove it we consider multigraphs and the configuration model of Bollobas \cite{Bol80, Bol01}.

Section~\ref{the typical structure of graphs with bounded maximum degree}
deals with the typical structure of members of $\mbG_{n,R}$ as $n \to \infty$.
The results are summarised in 
Theorem~\ref{theorem about typical structure}.
Intuitively speaking, there are three types of
``rare'' events, which are vertices with degree $R-2$,
short cycles and short paths with endpoints of degree $R-1$.
They are, asymptotically, Poisson distributed and independent of each other.
Moreover, these rare events typically occur far from each other.
For example, the distance between any vertex with degree $R-2$ and
any path of length 3 (say) with both endpoints of degree $R-1$
will almost surely be large.
Also, the ``semi-rare'' event
of being a vertex with degree $R-1$ almost surely occurs far from any rare event.
The proofs use the configuration model of Bollob\'{a}s~\cite{Bol80, Bol01} and the fact,
given by Theorem~\ref{theorem about typical degree distribution},
that we can restrict attention to degree sequences with certain properties.

In Section~\ref{First-order limit laws},
we use Theorems~\ref{theorem about typical degree distribution}
and~\ref{theorem about typical structure}
to prove a labelled logical limit law for $\mbG_{n,R}$.
In the cases $R = 0$ and $R = 1$ the result is known (and trivial),
in both the labelled and unlabelled context,
so we need only consider the case $R \geq 2$.
If $R \geq 5$ then we can also conclude that $\mbG_{n,R}$
satisfies an unlabelled limit law, by applying
a result of McKay and Wormald~\cite{MW84}, saying
that in both the labelled and unlabelled cases a random member of 
$\mbG_{n,R}$ almost surely has no non-trivial automorphism,
and a general model theoretic result concerning labelled and
unlabelled probabilities.
We cannot apply the main result of \cite{Lyn05} to get
the limit laws of this paper, because the degree sequences 
that are typical of members of $\mbG \in \mbG_{n,R}$ do
not satisfy condition~(1)(b) of \cite{Lyn05}.
Moreover, the structural results of 
Sections~\ref{the distribution of degrees}
and~\ref{the typical structure of graphs with bounded maximum degree}
are needed for understanding logical limit laws of $\mcH$-free graphs in \cite{Kop12}.

\medskip

\noindent
{\bf Preliminaries, terminology and notation.}
See for example \cite{EF, Lib, Spe} for definitions of {\em first-order logic},
By the {\em language of graphs} we mean the first-order formulas which can
be built from the {\em vocabulary} (also called {\em signature}) which consists
of a binary relation symbol $E$ (denoting the edge relation) and the
identity symbol `='. More generally, a finite relational vocabulary consists of
a list of relation symbols $R_1, \ldots, R_k$ (and we can assume that $R_1$ is `=')
each with an associated arity. A {\em structure} for the language associated with this
vocabulary is a tuple $\mcS = (S, R_1^{\mcS}, \ldots, R_1^{\mcS})$ such that $S$ is a
set and $R_i^{\mcS} \subseteq S^{r_i}$ where $r_i$ is the arity of $R_i$.
From this point of view a graph is a structure for the language of graphs, so it has the
form $\mcG = (V, E^{\mcG})$ where $(a,b) \in E^{\mcG}$ implies $(b,a) \in E^{\mcG}$
since we consider undirected graphs here. Alternatively we may, as in graph theory,
view $E^{\mcG}$ as the set of edges (2-subsets of $V$) of $\mcG$.
For a first-order sentence (a formula with no free variables) $\varphi$ in the language
of graphs and a graph $\mcG$, $\mcG \models \varphi$ is shorthand for 
``$\mcG$ satisfies $\varphi$''. See for example \cite{EF, Lib} for a formal definition of `$\models$'.

The degree of a vertex $v$ in a graph $\mcG$ is denoted $\deg_{\mcG}(v)$, and the notation
$v \sim_{\mcG} w$ means that vertices $v$ and $w$ are adjacent in $\mcG$
(which is equivalent with $\mcG \models E(v,w)$).
By $\dist_{\mcG}(v,w)$ we mean the distance from $v$ to $w$ in $\mcG$.
For functions $f$ and $g$ we say that $f$ is {\em asymptotic} with $g$, denoted $f \sim g$, 
if $\lim_{x\to\infty}f(x)/g(x) = 1$. 
For positive integers $n$ we sometimes use the notation $[n] = \{1, \ldots, n\}$.
To make notation simpler we sometimes omit explicit use of the floor function 
$\lfloor \ \rfloor$. For a set $S$, $|S|$ denotes its cardinality.
By $\mbbN$ we denote the set of non-negative integers.

\section{Distribution of degrees}
\label{the distribution of degrees}

\noindent
For positive integers $n$ and $R$, $\mbG_{n,R}$ denotes the set of graphs with vertices $1, \ldots, n$
in which the degree of every vertex is at most $R$.

\begin{theor}\label{theorem about typical degree distribution}
Suppose that $R \geq 2$.
\begin{itemize}
\item[(i)] \ The proportion of $\mcG \in \mbG_{n,R}$ which have no vertex with degree less than $R-2$
approaches 1 as $n \to \infty$.
\item[(ii)] \ For every $\varepsilon > 0$, 
the proportion of $\mcG \in \mbG_{n,R}$ which have between $\sqrt{(R - \varepsilon)n}$ and 
$\sqrt{(R + \varepsilon)n}$ vertices with degree $R-1$ approaches 1 as $n \to \infty$.
\item[(iii)] \ For every $k \in \mbbN$, the proportion of $\mcG \in \mbG_{n,R}$ which
have exactly $k$ vertices with degree $R-2$ approaches
$$\frac{(R-1)^k \ e^{-(R-1)}}{k!} \qquad \text{ as } \ n \to \infty.$$
In other words, the number of vertices with degree $R-2$ has, asymptotically, a Poisson distribution
with mean $R-1$.
\end{itemize}
\end{theor}

\noindent
If $3 \leq r \leq R$ and
$\bar{d} = (d_r, \ldots, d_R)$ is a sequence of  nonnegative integers summing to $n$,
then let $\mbG_{n,\bar{d}}$ be the set of all graphs $\mcG$ with vertices $1, \ldots, n$ such that,
for $i = r, \ldots, R$, $\mcG$ has exactly $d_i$ vertices with degree $i$.
Wormald (Theorem~1 in \cite{Wor81}) and McKay and Wormald 
(Corollaries~3.8 and~3.10 in \cite{MW84}) have proved the following:

\begin{theor}\label{result of Wormald} \cite{Wor81, MW84} 
If $3 \leq r \leq R$ and
$\bar{d} = (d_r, \ldots, d_R)$ is a sequence of  nonnegative integers summing to $n$,
then 
\begin{itemize}
\item[(i)] the number of $r$-connected $\mcG \in \mbG_{n,\bar{d}}$ is 
$|\mbG_{n,\bar{d}}| \big(1 - O\big(n^{2-r}\big)\big)$ where $O( \ )$ denotes a bound
depending only on $R$, and
\item[(ii)] the expected number of nontrivial automorphisms of a random $\mcG \in \mbG_{n,\bar{d}}$ 
is $O(n^{2-r})$ where $O( \ )$ denotes a bound depending only on $R$.
\end{itemize} 
Statement (ii) is true also if we count graphs in $\mbG_{n,\bar{d}}$ only up to
isomorphism, or in other words, if we consider unlabelled graphs with $n$ vertices and
degree sequence $\bar{d}$.
\end{theor}

\noindent
Let us call a graph {\em rigid} if it has only one automorphism, the trivial one.
Theorem 4.3.4 in~\cite{EF}, if stated in the context of graphs, says the following.

\begin{theor}\label{rigidity implies that labelled and unlabelled probabilities coincide}\cite{EF}
For every positive integer $n$, let $\mbH_n \subseteq \mbG_n$ be sets of graphs with 
vertices $1, \ldots, n$.
Suppose that, for every $n$, 
if $\mcG \in \mbH_n$, $\mcH$ has vertices $1, \ldots, n$ and is isomorphic to $\mcG$,
then $\mcH \in \mbH_n$; 
and if $\mcG \in \mbG_n$, $\mcH$ has vertices $1, \ldots, n$ and is isomorphic to $\mcG$,
then $\mcH \in \mbG_n$.
For every $n$ and $\mcG \in \mbG_n$, let $[\mcG]$ denote the equivalence class to which
$\mcG$ belongs with respect to the isomorphism relation on $\mbG_n$.
Let $\mb{RIG}(\mbG_n)$ be the set of $\mcG \in \mbG_n$ which are rigid.
\begin{align*}
\text{If } \ \
&\lim_{n\to\infty} \ \frac{|\{[\mcG] : \mcG \in \mb{RIG}(\mbG_n)\}|}{|\{ [\mcG] : \mcG \in \mbG_n\}|} \ = \ 1
\quad \text{ then } \\
&\lim_{n\to\infty}  \ \frac{|\{[\mcG] : \mcG \in \mbH_n\}|}{|\{[\mcG] : \mcG \in \mbG_n\}|} \ = \ 
\lim_{n\to\infty} \ \frac{|\mbH_n|}{|\mbG_n|},
\end{align*}
if any one of the last two limits exist.
\end{theor}

\noindent
We say that a graph has {\em connectivity} $k$ if it is $k$-connected but not $(k+1)$-connected.
(By $k$-connected we mean $k$-vertex connected.)
By combining Theorems~\ref{theorem about typical degree distribution},
~\ref{result of Wormald} 
and~\ref{rigidity implies that labelled and unlabelled probabilities coincide}
we get the following information about 
typical members of $\mbG_{n,R}$ for large $n$.

\begin{cor}\label{corollary to main theorem}
Suppose that $R \geq 5$.
The proportion of graphs $\mcG \in \mbG_{n,R}$ 
(or of isomorphism classes $[\mcG]$ where $\mcG \in \mbG_{n,R}$)
which have the following two properties approaches 1 as 
$n$ approaches infinity:
\begin{itemize}
\item[(i)] The connectivity of $\mcG$ is $R-2$ or $R-1$.
\item[(ii)] $\mcG$ is rigid.
\end{itemize}
\end{cor}

\noindent
Corollary~\ref{corollary to main theorem} will be used to derive an
unlabelled logical limit law in Section~\ref{First-order limit laws}.

\subsection{Multigraphs and configurations}\label{configurations and multigraphs}

Theorem~\ref{theorem about typical degree distribution} will be proved
via the use of multigraphs and so-called configurations (defined below).
The corresponding statements (of Theorem~\ref{theorem about typical degree distribution})
for multigraphs are given by Lemmas~\ref{no vertices with degree less than R-2 for multigraphs},
~\ref{distribution of vertices with degree R-1 for multigraphs}
and~\ref{Poisson distribution degree R-2 for multigraphs} below.
In Section~\ref{configurations and simple graphs} we explain how 
Theorem~\ref{theorem about typical degree distribution} follows from these lemmas
and from asymptotic estimates, of Bollobas \cite{Bol80, Bol01}, 
of probabilities of small cycles in configurations.
We fix an arbitrary integer $R \geq 1$.

\begin{defin}\label{definition of multigraph}{\rm
A {\em multigraph} is an object of the form $(V, f)$
where $f : U \to \{0,1,2, \ldots\}$ and 
$U = \{W \subseteq V : 1 \leq |W| \leq 2\}$.
If $\mcM = (V, f)$ is a multigraph and $v \in V$ then we define
$\displaystyle deg_\mcM(v) = 2f(\{v\}) + \sum_{\substack{w \in V \\ w \neq v}} f(\{v,w\})$ 
and call $deg_\mcM(v)$ the {\em degree of $v$} (in $\mcM$).
}\end{defin}

\noindent
Suppose that $(V,f)$ is a multigraph, $\{v,w\} \in U = \{W \subseteq V : 1 \leq |W| \leq 2\}$
and $f(\{v,w\}) = k$. Intuitively this means that there are exactly $k$ different edges between $v$ and $w$.
If $v = w$ then it means that there are exactly $k$ different loops which begin and end in $v$.
If $f(W) \leq 1$ for all $W \in U$ and $f(W) = 0$ whenever $|W| = 1$,
then $(V,f)$ corresponds to the graph $(V,E)$ where, for all distinct $v, w \in V$,
$\{v,w\} \in E$ if and only if $f(\{v,w\}) = 1$.

\begin{defin}\label{definition of set of multigraphs}{\rm
(i) Let $\mb{MG}_{n,R}$ be the set of multigraphs with vertices $1, \ldots, n$ such
that every vertex has degree at most $R$.\\
(ii) For every sequence $\bar{d} = (d_0, \ldots, d_R)$ of non-negative integers
such that $\sum_{i=0}^R d_i = n$ 
we let $\mb{MG}_{n, \bar{d}}$ denote the set of multigraphs with vertices
$1, \ldots, n$ such that there are exactly $d_i$ vertices with degree $i$ for $i = 0, \ldots, R$.\\
(iii) Let $N(d_0, \ldots, d_R) = |\mb{MG}_{n,\bar{d}}|$.
}\end{defin}

\begin{defin}\label{definition of configuration}{\rm
Let $E$ and $F$ be binary relation symbols and let $L_c$ be the first-order language
with vocabulary $\{=,E,F\}$.
A {\em configuration} is an $L_c$-structure $\mcC = (C, E^{\mcC}, F^{\mcC})$ such that
$F^{\mcC}$ is an equivalence relation on $C$ and $E^{\mcC}$ is an equivalence relation on $C$ such that
every $E^{\mcC}$-class has cardinality 2. 
The later condition implies that $|C|$ must be even (or infinite).
An $E^{\mcC}$-class of a configuration $\mcC$ will often be called an {\em edge} (of $\mcC$).
We say that $\mcC = (C, E^{\mcC}, F^{\mcC})$ is an {\em $n$-configuration} 
if $F^{\mcC}$ has exactly $n$ nonempty classes.
}\end{defin}

\noindent
Observe that if $\mcC = (C, E^{\mcC}, F^{\mcC})$ is a configuration, then 
$(C, E^{\mcC})$ is a graph which is a {\em complete matching} (also called perfect matching),
i.e. every $c \in C$ is adjacent to exactly one member of $C$.

\begin{defin}\label{notation for number of matchings}{\rm
If $m > 0$ is an integer then let $\mr{M}(2m)$ denote the number of complete
matchings, i.e. equivalence relations such that every class has cardinality 2,
on a set with cardinality $2m$.
}\end{defin}

\noindent
It is easy to see that
\begin{equation*}
\mr{M}(2m) = \binom{2m}{2}\binom{2m-2}{2} \cdots \binom{2}{2} \bigg/ m! \ = \ 
\frac{(2m)!}{m!2^m}.
\end{equation*}
By Stirlings approximation,
\begin{equation}\label{stirlings apporoximation of matchings}
\frac{(2m)!}{m! 2^m} \ = \ \bigg(\frac{2m}{e}\bigg)^m \big(\sqrt{2} + o(1)\big).
\end{equation}
Also note that for any fixed integer $p$ 
we have $m! / (m-p)! \sim m^p$ which gives
\begin{equation}\label{approximation of quotient of matchings}
\frac{\mr{M}(2m)}{\mr{M}(2m - 2p)} \ = \ 
\frac{(2m)!}{m! 2^m} \ \Bigg/ \ \frac{(2(m-p))!}{(m-p)! 2^{m-p}} 
\ \sim \ (2m)^p.
\end{equation}

\begin{defin}\label{definition of set of configurations}{\rm
We will always assume that $W_1, \ldots, W_n$ is a sequence of disjoint sets 
such that $\sum_{i=1}^n |W_i| = 2m$ for some integer $m > 0$.
By $\mbC(W_1, \ldots, W_n)$ we denote the set of configurations $\mcC = (C, E^{\mcC}, F^{\mcC})$
such that $C = \bigcup_{i=1}^n W_i$ and $W_1, \ldots, W_n$ enumerates all the $F^{\mcC}$-classes
(so if some set $W_i$ other than an $F^{\mcC}$-class appears in the sequence then $W_i$ is empty).
}\end{defin}

\noindent
Observe that if $2m = \sum_{i=1}^n |W_i|$, then 
\begin{equation}\label{number of configurations on Ws}
|\mbC(W_1, \ldots, W_n)| \ = \ \mr{M}(2m) \ = \ \frac{(2m)!}{m!2^m}.
\end{equation}

\begin{defin}\label{definition of multigraph image of a configuration}{\rm
Let $\mcC \in \mbC(W_1, \ldots, W_n)$.
The {\em multigraph image} of $\mcC$, denoted $\mr{Graph}(\mcC)$,
is the multigraph with vertices $1, \ldots, n$ 
such that, for all $i,j \in \{1, \ldots, n\}$, 
the number of edges between $i$ and $j$ is the same as the number of
$E^{\mcC}$-classes which have non-empty intersection with both $W_i$ and $W_j$.
($E^{\mcC}$-classes which are included in a single $W_i$, if they exist, give rise to loops
in $\mr{Graph}(\mcC)$.)
}\end{defin}

\noindent
Observe that for every $\mcC \in \mbC(W_1, \ldots, W_n)$ and every $i = 1, \ldots, n$,
the degree of $i$ in $\mr{Graph}(\mcC)$ is $|W_i|$.
Moreover, if $\mcG$ is a multigraph with vertices $1, \ldots, n$ such that, for $i = 1, \ldots,n$,
$i$ has degree $|W_i|$, then there are exactly $\prod_{i=1}^n|W_i|!$ configurations
$\mcC \in \mbC(W_1, \ldots, W_n)$ such that $\mr{Graph}(\mcC) = \mcG$, because
$\mr{Graph}(\mcC') = \mr{Graph}(\mcC)$ if and only if there is an isomorphism from
$\mcC'$ to $\mcC$ which preserves every $W_i$ setwise.
It follows that if $\bar{d} = (d_0, \ldots, d_R)$ is a sequence of non-negative integers
such that $\sum_{i=0}^R d_i = n$ and $\sum_{i=0}^R id_i = 2m$, then, 
using~(\ref{number of configurations on Ws}), we have
\begin{equation}\label{number of multigraphs with a given degree sequence}
N(d_0, \ldots, d_R) \ = \ 
\binom{n}{d_0, \ldots, d_R} \frac{(2m)!}{m!2^m} \Bigg/ \prod_{i=0}^R (i!)^{d_i}.
\end{equation}

{\em From now on assume that $\bar{d} = (d_0, \ldots, d_R)$ is a sequence of non-negative integers
such that $\sum_{i=0}^R d_i = n$ and $\sum_{i=0}^R id_i = 2m$ for some integer $m > 0$.}
We are interested in the
probability that a random $\mcG \in \mb{MG}_{n,R}$ (drawn uniformly) has this degree sequence. 
In other words, we consider the proportion
\begin{equation*}
p(d_0, \ldots, d_R) \ = \ \frac{N(d_0, \ldots, d_R)}{|\mb{MG}_{n,R}|} \ = \ 
\frac{N(d_0, \ldots, d_R)}{\sum N(d'_0, \ldots, d'_R)},
\end{equation*}
where the sum ranges over all sequences $(d'_0, \ldots, d'_R)$ of non-negative integers
such that $\sum_{i=0}^R d'_i = n$.

If $n$ is even and all vertices have degree $R$, that is, if $d_i = 0$ for $i = 0, \ldots, R-1$,
we get $2m = Rn$ and, by~(\ref{number of multigraphs with a given degree sequence})
and~(\ref{stirlings apporoximation of matchings}),
\begin{equation}\label{when all vertices have degree R}
N(0,\ldots, 0, n) \ = \ \bigg(\frac{Rn}{e}\bigg)^{Rn/2} \frac{1}{(R!)^n} \big(\sqrt{2} + o(1)\big) 
\ \geq \ C_1^n n^{Rn/2}
\end{equation}
for some constant $C_1 > 0$.
If $n$ is odd, then, with $d_0 = 1$, $d_i = 0$ for $i = 1, \ldots, R-1$ and $d_R = n-1$, we have
$2m = R(n-1)$ and get 
in a similar way
\begin{equation}\label{when all but one vertex has degree R}
N(1,0,\ldots,0, n-1) \ \geq \ C_2^n n^{Rn/2}
\end{equation}
for some constant $C_2 > 0$.

Let $\varepsilon > 0$.
We now estimate the proportion of $\mcG \in \mb{MG}_{n,R}$ for which $2m \leq (R-\varepsilon)n$.
Summing over all $(d_0, \ldots, d_R)$ for which $\sum_{i=0}^R id_i = 2m \leq (R-\varepsilon)n$
gives, for large enough $n$ and using~(\ref{when all vertices have degree R})
and~(\ref{when all but one vertex has degree R}),
\begin{align*}
&\sum_{2m \leq (R-\varepsilon)n} N(d_0, \ldots, d_R) \ \leq \
\sum_{2m \leq (R-\varepsilon)n} 
\binom{n}{d_0, \ldots, d_R} \frac{(2m)!}{m!2^m} \Bigg/ \prod_{i=0}^R (i!)^{d_i} \\
\leq \ 
&(R+1)^n \bigg(\frac{(R-\varepsilon)n}{e}\bigg)^{(R-\varepsilon)n/2} \cdot 2
\ = \ 
2(R+1)^n \bigg(\frac{(R-\varepsilon)n}{e}\bigg)^{Rn/2} 
\bigg(\frac{(R-\varepsilon)n}{e}\bigg)^{-\varepsilon n/2} \\
\ \leq \
&\Big(C_3 n^{-\varepsilon/2}\Big)^n N, \ \ 
\text{where $N = N(0,\ldots,0,n)$ or $N = N(1,0, \ldots, 0,n-1)$} \\
&\text{depending on whether
$n$ is even or odd and $C_3 > 0$ is a constant.}
\end{align*}

\noindent
It follows that
\begin{equation*}
\frac{\sum_{2m \leq (R-\varepsilon)n} N(d_0, \ldots, d_R)}{|\mb{MG}_{n,R}|} \ \leq \ 
\Big(C_3 n^{-\varepsilon / 2}\Big)^n \ \to \ 0 \quad \text{ as }  n \to \infty.
\end{equation*}
Therefore we will assume that $2m > (R-\varepsilon)n$. Since $\varepsilon > 0$ 
is arbitrary we may in fact assume that $2m > (R - o(1))n$.
Since
\begin{equation*}
Rn - 2m \ = \ Rn - \sum_{i=0}^R id_i \ = \ \sum_{i=0}^R (R-i)d_i
\end{equation*}
and $2m > (R - o(1))n$ now implies $Rn - 2m \leq o(n)$ we get
\begin{equation*}
\frac{\sum_{i=0}^R (R-i)d_i}{n} \ \leq \ o(1)
\end{equation*}
which in turn implies that
\begin{equation}\label{can assume that all d-i except d-R is o(n)}
\text{$d_i = o(n)$ for all $i = 0, \ldots, R-1$.}
\end{equation}

\noindent
Next, we aim at showing that in the typical case we have $d_i = 0$
for all $i = 0, \ldots, R-3$. So until further notice we assume that $R \geq 3$ 
since otherwise there is nothing to prove.
Fix any $i \leq R-3$ and suppose that $d_i \geq 1$.
Now let
\begin{align*}
&d'_i = d_i - 1,\\
&d'_{i+2} = d_{i+2} + 1 \ \text{ and}\\
&d'_j = d_j \ \text{ for all } j \notin \{i, i+2\}.
\end{align*}
With $2m' = \sum_{j=0}^R jd'_j$ we have 
$2m' = 2m + 2$ (where $2m = \sum_{j=0}^R jd_j$).
Now we get, by~(\ref{number of multigraphs with a given degree sequence})
and~(\ref{approximation of quotient of matchings}),
\begin{align*}
&\frac{N(d_0, \ldots, d_R)}{N(d'_0, \ldots, d'_R)} \\
= \ \ 
&\binom{n}{d_0, \ldots, d_R} \ \frac{(2m)!}{m!2^m} \ \bigg(\prod_{j=0}^R (j!)^{d_j}\bigg)^{-1} \
\Bigg/ \
\binom{n}{d'_0, \ldots, d'_R} \ \frac{(2m')!}{m'!2^{m'}} \ \bigg(\prod_{j=0}^R (j!)^{d'_j}\bigg)^{-1} \\
\sim \ \ 
&\frac{d'_i! \ d'_{i+2}!}{d_i! \ d_{i+2}!} \ \frac{1}{2m} \ \prod_{j=0}^R (j!)^{d'_j - d_j}
\ = \ 
\frac{d_{i+2} + 1}{d_i} \ \frac{1}{2m} \ (i+2)(i+1).
\end{align*}
As we assume that $2m > (R - o(1))n$ we may assume that $2m \geq n$.
If $i \leq R-3$ then, by~(\ref{can assume that all d-i except d-R is o(n)}),
$d_{i+2} = o(n)$ so we may assume that $d_{i+2} + 1 \leq f(n)$ where $f(n) = o(n)$,
which gives
\begin{equation*}
\frac{N(d_0, \ldots, d_R)}{N(d'_0, \ldots, d'_R)} \ \leq \ 
\frac{f(n) R^2}{n}.
\end{equation*}
Hence, for $i \leq R-3$,
\begin{equation*}
p(d_0, \ldots, d_R) 
\ \leq \ 
\frac{R^2 f(n)}{n} p(d'_0, \ldots, d'_R),
\end{equation*}
and if we sum over all $(d_0, \ldots, d_R)$ for which 
$d_i \geq 1$, $d_{i+2} \leq f(n)$ and $2m \geq n$, we get
\begin{align*}
\sum p(d_0, \ldots, d_R) \ 
&\leq \ 
\frac{R^2 f(n)}{n} \sum p(d'_0, \ldots, d'_R) \\ 
&\leq \ 
\frac{R^2 f(n)}{n} \ = \ o(1) \quad \text{(since $f(n) = o(n)$)}.
\end{align*}
This implies that the proportion of $\mcG \in \mb{MG}_{n,R}$ which have a vertex
with degree less than $R-2$ approaches 0 as $n \to \infty$.
Therefore we assume, from now on, that
\begin{equation}\label{d-i is zero for all i < R - 2}
\text{$d_i = 0$ for all $i = 0, \ldots, R-3$.}
\end{equation}
We also state this conclusion as:

\begin{lem}\label{no vertices with degree less than R-2 for multigraphs}
The proportion of $\mcG \in \mb{MG}_{n,R}$ such that $\mcG$ has a vertex with
degree less than $R-2$ approaches 0 as $n \to \infty$.
\end{lem}

The next step is to find a lower upper bound
(than in~(\ref{can assume that all d-i except d-R is o(n)}))
for $d_{R-2}$ in the typical case, so we assume that $R \geq 2$ until further notice.
Suppose that $g(n) \to \infty$ as $n \to \infty$ and
suppose that $d_{R-2} \geq g(n)$.
We argue as before, with $i = R-2$, so
$d'_{R-2} = d_{R-2} - 1$, $d'_R = d_R + 1$ and $d'_j = d_j$ if $j \notin \{R-2, R\}$.
Just as before, and as $2m > (R - o(1))n$ where $R \geq 2$, we get
\begin{equation*}
\frac{N(d_0, \ldots, d_R)}{N(d'_0, \ldots, d'_R)} \ \sim \ 
\frac{d_R + 1}{d_{R-2}} \ \frac{1}{2m} \ R(R-1) \ < \ 
\frac{R^2}{d_{R-2}} \ \leq \ \frac{R^2}{g(n)}.
\end{equation*}
Hence $p(d_0, \ldots, d_R) \leq R^2 p(d'_0, \ldots, d'_R) / g(n)$ and if we sum over
all $(d_0, \ldots, d_R)$ such that $d_{R-2} \geq g(n)$ 
we get 
\begin{equation*}
\sum p(d_0, \ldots, d_R) \ \leq \ \frac{R^2}{g(n)}.
\end{equation*}
Therefore we may assume that 
\begin{equation}\label{d-R-2 cannot grow}
\text{$d_{R-2} \leq g(n)$ for every $g$ such that $\lim_{n\to\infty} g(n) = \infty$.}
\end{equation}

The next step is to show that in the typical case $d_{R-1} = O(\sqrt{n})$.
Assume that $R \geq 1$ and that $d_{R-1} \geq 2$.
Let 
\begin{align*}
&d'_{R-1} = d_{R-1} - 2, \\
&d'_R = d_R + 2 \ \text{ and}\\
&d'_j = d_j \text{ if } j \notin \{R-1, R\}.
\end{align*}
Then $2m' = 2m + 2$ where $2m' = \sum_{j=0}^R id'_j$ and
\begin{align*}
\frac{N(d_0, \ldots, d_R)}{N(d'_0, \ldots, d'_R)} \ 
&\sim \ 
\frac{(d_R + 1)(d_R + 2)}{d_{R-1}(d_{R-1} - 1)} \ \frac{1}{2m} \ R^2 \\
&\leq \ 
\frac{(d_R + 1)n}{d_{R-1}(d_{R-1} - 1)} \ \frac{1}{Rd_R} \ R^2 \ 
\leq \ 
\frac{4Rn}{d_{R-2}^2}.
\end{align*}
It follows that 
$$p(d_0, \ldots, d_R) \ \leq \ \frac{4Rn}{d_{R-1}^2} \ p(d'_0, \ldots, d'_R).$$
Let $\mbbP(d_{R-1} \geq k)$ denote the proportion of $\mcG \in \mb{MG}_{n,R}$
which have at least $k$ vertices with degree $R-1$. Then for every $k \geq 2$, by summation over
all $(d_0, \ldots, d_R)$ such that $d_{R-1} \geq k$, we get 
\begin{equation}\label{reduction from k to k-2}
\mbbP(d_{R-1} \geq k) \ \leq \ \frac{4Rn}{k^2} \ \mbbP(d_{R-1} \geq k-2).
\end{equation}
If $C > 2\sqrt{R}$ then, by iterating~(\ref{reduction from k to k-2}) $\sqrt{n}$ times,
we get that  the proportion of $\mcG \in \mb{MG}_{n,R}$ which have at least 
$(C + 2)\lceil\sqrt{n}\rceil$ vertices with degree $R-1$ is
\begin{equation*}
\mbbP\Big(d_{R-1} \geq (C + 2)\lceil\sqrt{n}\rceil\Big) \ \leq \ 
\Bigg(\frac{4Rn}{\big(C\sqrt{n}\big)^2}\Bigg)^{\sqrt{n}} \ = \ \bigg(\frac{4R}{C^2}\bigg)^{\sqrt{n}}
\ \to \ 0 \quad \text{ as } n \to \infty.
\end{equation*}
Therefore we may assume, in addition to the previous 
assumptions~(\ref{d-i is zero for all i < R - 2}) and~(\ref{d-R-2 cannot grow}), that
\begin{equation}\label{number of vertices with degree R-1 is Ordo of square rot n}
d_{R-1} = O(\sqrt{n}).
\end{equation}

\noindent
We use the notation $(n)_k = n(n-1)\cdots (n - k + 1)$. 
By equation~(3.8) in~\cite{Ben74}, for example, we have:
\begin{equation}\label{approximation of falling factorial}
\text{If $k = O(\sqrt{n})$, then } \ (n)_k = n^k \exp\bigg( -\frac{k^2}{2n} + o(1)\bigg).
\end{equation}
To make the argument here self contained we prove~(\ref{approximation of falling factorial}).
If $k = O(\sqrt{n})$ then
\begin{align*}
\ln\bigg(\frac{(n)_k}{n^k}\bigg) \ 
&= \  
\ln\frac{n!}{(n-k)!n^k} \ = \ 
\ln\Bigg( \prod_{i=0}^{k-1} \frac{n-i}{n} \Bigg) \\ 
&= \ 
\ln\Bigg( \prod_{i=0}^{k-1} \Big(1 - \frac{i}{n} \Big)\Bigg) \ = \ 
\sum_{i=0}^{k-1} \ln\Big(1 - \frac{i}{n}\Big),
\end{align*}
and for large $n$ and $i < k = O(\sqrt{n})$, $i/n$ is close to 0,
so by Taylor's approximation we get
\begin{align*}
\sum_{i=0}^{k-1} \ln\Big(1 - \frac{i}{n}\Big) \ &= \ 
\sum_{i=0}^{k-1}\bigg( - \frac{i}{n} \ + \ O\bigg( \Big(\frac{i}{n}\Big)^2 \bigg) \bigg) \\ 
&= \ 
-\frac{k(k-1)}{2n} \ + \ O\bigg(\frac{k^3}{n^2}\bigg) \ = \ 
-\frac{k^2}{2n} \ + \ \frac{k}{2n} \ + \ O\bigg(\frac{k^3}{n^2}\bigg) \\ 
&= \ 
-\frac{k^2}{2n} \ + \ o(1), \quad \text{ since } k = O(\sqrt{n}).
\end{align*}
Hence, $\ln\Big(\frac{(n)_k}{n^k}\Big) = -\frac{k^2}{2n} \ + \ o(1)$
which immediately gives~(\ref{approximation of falling factorial}).

Recall the exact formula~(\ref{number of multigraphs with a given degree sequence})
for $N(d_0, \ldots, d_R)$ and remember the 
assumptions~(\ref{d-i is zero for all i < R - 2}),
~(\ref{d-R-2 cannot grow})
and~(\ref{number of vertices with degree R-1 is Ordo of square rot n}).
Hence $d_i = 0$ for $i = 0, \ldots, R-3$, $d_{R-2} = o(\sqrt{n})$,
$d_{R-1} = O(\sqrt{n})$, $d_R = n - d_{R-1} - d_{R-2}$ and
$2m = Rn - d_{R-1} - 2d_{R-2}$ (since, by definition, $2m = \sum_{i=0}^R id_i$).
By also using~(\ref{approximation of falling factorial}) we get 
\begin{align}\label{simplified formula for number of multigraphs}
&N(d_0, \ldots, d_R) \ = \\
&= \ 
\binom{n}{d_{R-2}, d_{R-1}, d_R} \cdot
\frac{(2m)!}{m! 2^m} \cdot 
\frac{1}{((R-2)!)^{d_{R-2}} \ ((R-1)!)^{d_{R-1}} \ (R!)^{d_R}} \nonumber \\ 
&= \ 
\frac{(n)_{d_{R-2} + d_{R-1}}}{d_{R-2}! \ d_{R-1}!} \cdot
\frac{(2m)!}{m! 2^m} \cdot
\frac{(R(R-1))^{d_{R-2}} \ R^{d_{R-1}}}{(R!)^n} \nonumber \\
&= \ 
\frac{n^{d_{R-2} + d_{R-1}} \ \exp\Big( - \frac{(d_{R-2} + d_{R-1})^2}{2n} \ + \ o(1)\Big)}
{d_{R-2}! \ d_{R-1}!} \cdot
\frac{(2m)!}{m! 2^m} \cdot
\frac{(R(R-1))^{d_{R-2}} \ R^{d_{R-1}}}{(R!)^n}. \nonumber
\end{align}

\noindent
From $2m = Rn - d_{R-1} - 2d_{R-2}$ and~(\ref{approximation of falling factorial})
it follows that
\begin{equation}\label{semifactorial of 2m}
\frac{(2m)!}{m! 2^m} \ = \ 
\frac{(Rn)!}{\lfloor Rn/2\rfloor! 2^{Rn/2}} \cdot
\frac{\exp\Big( \frac{(d_{R-1} + 2d_{R-2})^2}{4Rn} \ + \ o(1)\Big)}{(Rn)^{(d_{R-1} + 2d_{R-2})/2}}.
\end{equation}
By combining~(\ref{simplified formula for number of multigraphs})
and~(\ref{semifactorial of 2m}) we get
\begin{align}\label{final asymptotic formula for number of multigraphs}
N(d_0, \ldots, d_R) \ = \ 
&C_n \cdot \frac{(R-1)^{d_{R-2}}}{d_{R-2}!} \cdot
\frac{\big(\sqrt{Rn}\big)^{d_{R-1}}}{d_{R-1}!} \cdot \\
&\exp\Bigg( - \frac{(d_{R-1} + d_{R-2})^2}{2n} \ + \ 
\frac{(d_{R-1} + 2d_{R-2})^2}{4Rn} \ + \ o(1) \Bigg), \nonumber
\end{align}
where
$$C_n \ = \ \frac{(Rn)!}{\lfloor Rn/2\rfloor! 2^{Rn/2}} \cdot
\frac{1}{(R!)^n}$$
depends only on $n$.
Note that since $d_{R-2} = o(\sqrt{n})$ and $d_{R-1} = O(\sqrt{n})$ it follows that the
expression `$\exp( \ldots )$' in~(\ref{final asymptotic formula for number of multigraphs})
is bounded as $n\to\infty$.

\begin{lem}\label{distribution of vertices with degree R-1 for multigraphs}
(i) For every $c < R$ the proportion of $\mcG \in \mb{MG}_{n,R}$ with less than $\sqrt{cn}$
vertices with degree $R-1$ approaches 0 as $n\to\infty$.\\
(ii) For every $c > R$ the proportion of $\mcG \in \mb{MG}_{n,R}$ with more than $\sqrt{cn}$
vertices with degree $R-1$ approaches 0 as $n\to\infty$.
\end{lem}

\noindent
{\bf Proof.}
(i) Let $c < R$. Since there are not more than $n$ possibilites for $d_{R-1}$,
it suffices to prove that for all large enough $n$,
$$\frac{N(0,\ldots,0, d_{R-2}, \lfloor \sqrt{cn} \rfloor, d_R)}
{N(0,\ldots,0, d_{R-2}, \lfloor \sqrt{cn} \rfloor + 2\lfloor n^{1/4} \rfloor, d_R)} \ = \ 
\Bigg( \sqrt{\frac{c}{R}} + o(1) \Bigg)^{2\lfloor n^{1/4} \rfloor} 
\cdot O(1).$$
We may assume that $d_{R-2} = o(\sqrt{n})$ so 
by~(\ref{final asymptotic formula for number of multigraphs}) 
the above quotient equals
\begin{align*}
&(\sqrt{Rn})^{-2\lfloor n^{1/4} \rfloor} \cdot 
\frac{(\lfloor\sqrt{cn}\rfloor + 2\lfloor n^{1/4} \rfloor)!}{\lfloor \sqrt{cn} \rfloor !} \cdot
O(1) \ = \ 
(\sqrt{Rn})^{-2\lfloor n^{1/4} \rfloor} \cdot 
(\lfloor\sqrt{cn}\rfloor + 2\lfloor n^{1/4} \rfloor)_{2\lfloor n^{1/4} \rfloor} \\
\ \overset{(\ref{approximation of falling factorial})}{=} \
&(\sqrt{Rn})^{-2\lfloor n^{1/4} \rfloor} \ 
(\lfloor\sqrt{cn}\rfloor + 2\lfloor n^{1/4} \rfloor)^{2\lfloor n^{1/4} \rfloor} 
\cdot O(1)
\ = \ 
\Bigg(\frac{\lfloor\sqrt{cn}\rfloor + 2\lfloor n^{1/4} \rfloor}{\sqrt{Rn}}\Bigg)^{2\lfloor n^{1/4} \rfloor} 
\cdot O(1) \\
\ = \ 
&\Bigg( \sqrt{\frac{c}{R}} + o(1) \Bigg)^{2\lfloor n^{1/4} \rfloor} 
\cdot O(1).
\end{align*}
Part~(ii), where $c > R$, is proved in a similar way by considering
the quotient
$$\frac{N(0,\ldots,0, d_{R-2}, \lfloor \sqrt{cn} \rfloor, d_R)}
{N(0,\ldots,0, d_{R-2}, \lfloor \sqrt{cn} \rfloor - 2\lfloor n^{1/4} \rfloor, d_R)}. 
\qquad \qquad \square$$

\begin{lem}\label{Poisson distribution degree R-2 for multigraphs}
For every fixed $d_{R-2}$, the proportion of $\mcG \in \mb{MG}_{n,R}$ which
have exactly $d_{R-2}$ vertices with degree $R-2$ approaches 
$$\frac{(R-1)^{d_{R-2}} \ e^{-(R-1)}}{d_{R-2}!}.$$
\end{lem}

\noindent
{\bf Proof.}
Let $d_{R-2}$ be fixed and let $c > R$.
We know, by~(\ref{d-R-2 cannot grow})
and Lemma~\ref{distribution of vertices with degree R-1 for multigraphs}, 
that the proportion of $\mcG \in \mb{MG}_{n,R}$
such that $\mcG$ has more than $\lfloor n^{1/4} \rfloor$ vertices with degree $R-2$
or more than $\sqrt{cn}$ vertices with degree $R-1$,
approaches 0 as $n \to \infty$.
Therefore it suffices to prove that whenever $d_{R-1} \leq \sqrt{cn}$,
$$\frac{N(0,\ldots,0,d_{R-2},d_{R-1},d_R)}
{\sum_{k=0}^{\lfloor n^{1/4} \rfloor} N(0,\ldots,0,k,d_{R-1},d_R)} \ = \ 
(1 \pm o(1)) \frac{(R-1)^{d_{R-2}} \ e^{-(R-1)}}{d_{R-2}!}.$$
By~(\ref{final asymptotic formula for number of multigraphs}), the above quotient equals
\begin{align*}
(1 \pm o(1)) \frac{\frac{(R-1)^{d_{R-2}}}{d_{R-2}!}}
{\sum_{k=0}^{\lfloor n^{1/4} \rfloor} \frac{(R-1)^k}{k!}} \ = \ 
(1 \pm o(1)) \frac{\frac{(R-1)^{d_{R-2}}}{d_{R-2}!}}{e^{R-1}} \ = \ 
(1 \pm o(1)) \frac{(R-1)^{d_{R-2}} \ e^{-(R-1)}}{d_{R-2}!}.
\end{align*}
\hfill $\square$

\subsection{Configurations and graphs}\label{configurations and simple graphs}

In this section we explain why Lemmas~\ref{no vertices with degree less than R-2 for multigraphs},
\ref{distribution of vertices with degree R-1 for multigraphs}
and~\ref{Poisson distribution degree R-2 for multigraphs}
also apply to graphs without multiple edges or loops.

\begin{defin}\label{definition of path and cycle of a configuration}{\rm
Let $\mcC = (C, E^{\mcC}, F^{\mcC})$ be a configuration and let $p \geq 1$ be an integer.
By a {\em path (of $\mcC$) from $a \in C$ to $b \in C$ of length $p$}, or {\em $p$-path} from $a$ to $b$,
we mean sequence of $p$ distinct edges (i.e. $E^{\mcC}$-classes) $e_1, \ldots, e_p$
such that there are $F^{\mcC}$-classes 
$W_1, \ldots, W_{p+1}$ such $a \in e_1 \cap W_1$, $b \in e_p \cap W_{p+1}$,
$W_i \neq W_j$ if $i \neq j$ and $i \neq 1$ or if $i \neq j$ and $j \neq p+1$ 
(but we allow that $W_1 = W_{p+1}$)
and, for $i = 1, \ldots, p$, $e_i$ contains a member of $W_i$ and a member of $W_{i+1}$.
In the described situation we call $a$ and $b$ the {\em endpoints} of the path.
By a {\em $p$-cycle} (of $\mcC$) we mean a $p$-path from $a$ to $b$ for some $a$ and $b$
such that $(a,b) \in F^{\mcC}$; or in other words, a $p$-path from $a$ to $b$ where
$a,b \in W$ for some $F^{\mcC}$-class $W$.
}\end{defin}

\noindent
Recall Definition~\ref{definition of set of configurations}
of $\mbC(W_1, \ldots, W_n)$ and
Definition~\ref{definition of multigraph image of a configuration}
about the multigraph image $\mr{Graph}(\mcC)$ of a configuraion $\mcC$.
It is easy to see that that if $\mcC \in \mbC(W_1, \ldots, W_n)$ 
and $\mcC$ has no $p$-cycle for $p = 1,2$, then $\mr{Graph}(\mcC)$ has no loops
or multiple edges.

\begin{defin}\label{definition of set of configurations without 1,2-cycles}{\rm
By $\mbC'(W_1, \ldots, W_n)$ we denote the set of configurations \\
$\mcC \in \mbC(W_1, \ldots, W_n)$ which have no $p$-cycle for $p = 1,2$.
}\end{defin}

\noindent
The proof of Theorem~1 in~\cite{Bol80} 
(or Theorem~2.16 in~\cite{Bol01}) shows the following:

\begin{fact}\label{proportion of configurations whose images are graphs}{\rm \cite{Bol80, Bol01}
Let $W_1, \ldots, W_n$ be disjoint sets and for 
$i = 1, \ldots, R$, let 
$$d_i = |\{j : |W_j| = i\}|.$$
Suppose that for all $i < R-2$, $d_i = 0$ and $d_{R-2} + d_{R-1} = O(\sqrt{n})$.\\
(i) Then
$$\Bigg| \frac{|\mbC'(W_1, \ldots, W_n)|}{|\mbC(W_1, \ldots, W_n)|} \ - \
\exp\bigg(-\frac{R-1}{2} - \frac{(R-1)^2}{4}\bigg) \Bigg| \ = \ o(1) \quad \text{ as } n \to \infty,$$
where the bound `$o( \ )$' depends only on $R$.\\
(ii) From~(\ref{number of configurations on Ws})
and~(\ref{number of multigraphs with a given degree sequence}) 
it follows that if we let $\bar{d} = (d_0, \ldots, d_R)$ (with $d_i$ as above for $i = 1, \ldots, R$
and $d_0 = 0$), then
$$\Bigg| \frac{|\mbG_{n,\bar{d}}|}{|\mb{MG}_{n,\bar{d}}|} \ - \
\exp\bigg(-\frac{R-1}{2} - \frac{(R-1)^2}{4}\bigg) \Bigg| \ = \ o(1) \quad \text{ as } n \to \infty,$$
where the bound `$o( \ )$' depends only on $R$.\\
}\end{fact}

\noindent
From part~(ii) of the fact it immediately follows that 
Lemmas~\ref{no vertices with degree less than R-2 for multigraphs}
and~\ref{distribution of vertices with degree R-1 for multigraphs}
apply to graphs as well, that is, they remain true if we replace $\mb{MG}_{n,R}$
by $\mbG_{n,R}$.
To see why Lemma~\ref{Poisson distribution degree R-2 for multigraphs}
also holds with $\mbG_{n,R}$ in the place of $\mb{MG}_{n,R}$ consider its proof.
There we showed that 
$$\frac{N(0,\ldots,0,d_{R-2},d_{R-1},d_R)}
{\sum_{k=0}^{\lfloor n^{1/4} \rfloor} N(0,\ldots,0,k,d_{R-1},d_R)} \ = \ 
(1 \pm o(1)) \frac{(R-1)^{d_{R-2}} \ e^{-(R-1)}}{d_{R-2}!}.$$
By Fact~\ref{proportion of configurations whose images are graphs}~(ii),
if $\bar{d} = (0,\ldots,0,d_{R-2},d_{R-1},d_R)$, $d_{R-2} + d_{R-1} = O(\sqrt{n})$, then
$$|\mbG_{n,\bar{d}}| \ = \ (1 \pm o(1))
\exp\bigg(-\frac{R-1}{2} - \frac{(R-1)^2}{4}\bigg)
N(0,\ldots,0,d_{R-2},d_{R-1},d_R).$$
In the same way, for each $\bar{d}' = (0,\ldots,0,k,d_{R-1},d_R)$,
we also get (assuming $k + d_{R-1} = O(\sqrt{n})$)
$$|\mbG_{n,\bar{d}'}| \ = \ (1 \pm o(1))
\exp\bigg(-\frac{R-1}{2} - \frac{(R-1)^2}{4}\bigg)
N(0,\ldots,0,k,d_{R-1},d_R - k).$$
Therefore, the proof of Lemma~\ref{Poisson distribution degree R-2 for multigraphs}
implies that, for any fixed $d_{R-2}$, the proportion of $\mcG \in \mbG_{n,R}$
such that $\mcG$ has exactly $d_{R-2}$ vertices with degree $R-2$ approaches \\
$(R-1)^{d_{R-2}} e^{-(R-1)} \big/ d_{R-2}!$ as $n \to \infty$.
In other words, Lemma~\ref{Poisson distribution degree R-2 for multigraphs}
holds if $\mb{MG}_{n,R}$ is replaced by $\mbG_{n,R}$.
This completes the proof of 
Theorem~\ref{theorem about typical degree distribution}.

\section{The typical asymptotic structure}
\label{the typical structure of graphs with bounded maximum degree}

\noindent
In this section we study the typical structure of large graphs with maximum degree
$R$, where $R \geq 2$ is a fixed integer. The following theorem summarizes the results
that will be proved.

\begin{theor}\label{theorem about typical structure}
Let $R \geq 2$ be an integer.\\
(i) For every integer $k > 0$ the proportion of graphs $\mcG \in \mbG_{n,R}$ which have properties 
(1)--(4) below approaches 1 as $n \to \infty$:
\begin{enumerate}
\item If $p \leq k$ then $\mcG$ has no subgraph $\mcH$ with exactly $p$ vertices and more than $p$ edges.
It follows that whenever $p_1, p_2 \geq 3$ and $p_1 + p_2 + p_3 \leq k$, then $\mcG$ does not have
a $p_1$-cycle and a $p_2$-cycle such that there is a $p_3$-path from a vertex in
the first cycle to a vertex in the second cycle. 
\item If $p, q \leq k$ then there is no vertex $v$ with degree less than $R$ and $p$-path 
from $v$ to a vertex that belongs to a $q$-cycle. In particular, no $q$-cycle contains
a vertex of degree less than $R$.
\item There do not exist distinct vertices $v_1, v_2, v_3$ all of
which have degree at most $R-1$ such that for all distinct $i,j \in \{1,2,3\}$ there is a
path of length at most $k$ from $v_i$ to $v_j$.
\item There do not exist distinct vertices $v$ and $w$ such that $\deg_{\mcG}(v) \leq R-1$,
$\deg_{\mcG}(w) \leq R-2$ and there is a path of length at most $k$ from $v$ to $w$.
\end{enumerate}
(ii) Let $k \geq 3$ be an integer. 
There are positive $\lambda_3, \ldots, \lambda_k, \mu_1, \ldots, \mu_k \in \mbbQ$ 
such that for all $r_3, \ldots, r_k$, $s_1, \ldots, s_k \in \mbbN$ 
the proportion of $\mcG \in \mbG_{n,R}$ which, for $p = 3, \ldots, k$, have exactly 
$r_p$ $p$-cycles and, for $p = 1, \ldots, k$, have exactly $s_p$ $p$-paths with both endpoints
of degree $R-1$, approaches
$$\Bigg(\prod_{p = 3}^k \frac{(\lambda_p)^{r_p} \ e^{-\lambda_p}}{r_p!}\Bigg)
\Bigg(\prod_{p = 1}^k \frac{(\mu_p)^{s_p} \ e^{-\mu_p}}{s_p!}\Bigg)
\quad \text{ as } \ n \to \infty.$$
Moreover, the mentioned proportion approaches this value independently
of the number of vertices with degree $R-2$.\\
(iii) If $R \geq 3$ then for every $k \in \mbbN$ the proportion of $\mcG \in \mbG_{n,R}$
such that $\mcG$ has a connected component with at most $k$ vertices approaches 0 as $n \to \infty$.
\end{theor}

\noindent
Note that part~(ii) of the theorem states that if 
$X_p$ is the number of $p$-cycles, for $p = 3, \ldots, k$, and
$Y_p$ is the number of $p$-paths with endpoints of degree $R-1$, for $p = 1, \ldots, k$,
then the random variables $X_3, \ldots, X_p, Y_1, \ldots, Y_p$ are, asymptotically,
independent Poisson variables with means $\lambda_3, \ldots, \lambda_k, \mu_1, \ldots, \mu_k$, respectively,
which are described in detail by Lemma~\ref{distribution of pairs of vertices of degree R-1 on a given distance}.
If one omits the consideration of short paths with endpoints of degree $R-1$ in part~(ii),
then the resulting statement is a straightforward consequence of 
Theorem~\ref{theorem about typical degree distribution} 
and either one of Theorem~1 in~\cite{Bol80}, (the proof of) Theorem~2.16
in~\cite{Bol01}, or Corollary~1 in~\cite{Wor81b}. 
However, knowing the distribution of short paths with endpoints of degree $R-1$
is necessary in the proof of the limit law in Section~\ref{First-order limit laws}.

\subsection{Graphs with a given degree sequence}\label{graphs with a given degree sequence}

Fix an integer $R \geq 2$.
By Theorem~\ref{theorem about typical degree distribution}, 
\begin{itemize}
\item[\textbullet] \ the proportion of $\mcG \in \mbG_{n,R}$ without 
vertices with degree less than $R-2$ approaches 1 as $n \to \infty$,
\item[\textbullet] \ the number of vertices of $\mcG \in \mbG_{n,R}$ 
with degree $R-2$ has a Poisson distribution, asymptotically, and
\item[\textbullet] \ for every $\varepsilon > 0$, 
the proportion of $\mcG \in \mbG_{n,R}$ which have between $\sqrt{(R - \varepsilon)n}$ and 
$\sqrt{(R + \varepsilon)n}$ vertices with degree $R-1$ approaches 1 as $n \to \infty$.
\end{itemize}
Therefore, there are $\varepsilon_n$ and $\delta_n$ such that 
$\lim_{n\to\infty} \varepsilon_n = \lim_{n\to\infty} \delta_n = 0$,
the proportion of $\mcG \in \mbG_{n,R}$ without vertices of degree 0,
with at most $n^{1/4}$ vertices with degree $R-2$ and 
with between $\sqrt{(R - \delta_n)n}$
and $\sqrt{(R + \delta_n)n}$ 
vertices with degree $R-1$ is at least $1 - \varepsilon_n$.
Let $\mbG'_{n,R}$ be the set of all $\mcG \in \mbG_{n,R}$ such that 
$\mcG$ has no vertices of degree less than $R-2$, 
at most $n^{1/4}$ vertices of degree $R-2$ and
between $\sqrt{(R - \delta_n)n}$ and $\sqrt{(R + \delta_n)n}$ vertices with degree $R-1$.
It follows that $|\mbG'_{n,R}| \big/ |\mbG_{n,R}| \to 1$ as $n \to \infty$.

Since the number of edges of $\mcG = (V, E^{\mcG})$
equals $\frac{1}{2}\sum_{v \in V} \deg_{\mcG}(v)$ it follows that $\sum_{v \in V}\deg_{\mcG}(v)$ 
must be even. 
For every positive $n \in \mbbN$ and sequence of integers $\bar{d} = (d_1, d_2, \ldots, d_n)$
such that $0 \leq d_i \leq R$ for all $i = 1, \ldots, n$ and $\sum_{i=1}^n d_i$ is even, let 
$\mbG_{n,R}(\bar{d})$ be the set of graphs $\mcG$ with vertices $1, \ldots, n$ such that
$\deg_{\mcG}(i) = d_i$ for $i = 1, \ldots, n$. ($\mbG_{n,R}(\bar{d})$ is different
from $\mbG_{n,\bar{d}}$ in Section~\ref{the distribution of degrees}.)
A priori we do not know if $\mbG_{n,R}(\bar{d})$ is non-empty for every $\bar{d}$ such
that $\sum_{i=1}^n d_i$ is even. 
But for such degree sequenes $\bar{d}$ as we will consider 
(satisfying~(\ref{conditions on degree sequence}) below) it follows from 
Fact~\ref{proportion of configurations whose images are graphs} above that 
$\mbG_{n,R}(\bar{d})$ is indeed non-empty if, in addition, $\sum_{i=1}^n d_i$ is even and $n$ is sufficiently large.
Bollobas \cite{Bol80, Bol01} gives asymptotic estimates of $|\mbG_{n,R}(\bar{d})|$.

Now suppose that $\sigma_n$ are positive numbers such that $\lim_{n\to\infty}\sigma_n = 0$,
that $P$ is a property of graphs and that $0 < c < 1$.
Moreover, assume that whenever $n$ is large enough and 
$\bar{d} = (d_1, d_2, \ldots, d_n)$ satisfies that
\begin{align}\label{conditions on degree sequence}
&R - 2 \leq d_i \leq R \ \text{ for all } \ i = 1, \ldots, n, \ 
|\{i : d_i = R-2\}| \ \leq \ n^{1/4} \ \text{ and } \\ 
&\sqrt{(R - \delta_n)n} \ \leq \ |\{ i : d_i = R-1 \}| \ \leq \ \sqrt{(R + \delta_n)n} \nonumber
\end{align}
and $\mbG_{n,R}(\bar{d}) \neq \es$, then the proportion of 
$\mcG \in \mbG_{n,R}(\bar{d})$ which have property $P$ is at least $c - \sigma_n$
and at most $c + \sigma_n$.
If $\mbP_n$ is the set of all $\mcG \in \mbG'_{n,R}$ which have property $P$
and, for every $\bar{d}$ satisfying~(\ref{conditions on degree sequence}), 
$\mbP_n(\bar{d})$ is the set of $\mcG \in \mbG_{n,R}(\bar{d})$ which have property $P$,
then we get
\begin{align*}
&|\mbP_n| \ = \ \sum_{\bar{d}}|\mbP_n(\bar{d})| \ \geq \ 
(c - \sigma_n) \sum_{\bar{d}} |\mbG_{n,R}(\bar{d})| \ = \ 
(c - \sigma_n) |\mbG'_{n,R}| \\
&\text{and in a similar way } \ 
|\mbP_n| \ \leq \ (c + \sigma_n) |\mbG'_{n,R}|,
\end{align*}
where the sums range over all $\bar{d}$ which satisfy~(\ref{conditions on degree sequence}),
so $c - \sigma_n \leq |\mbP_n|\big/|\mbG'_{n,R}| \leq c + \sigma_n$.
It follows that the proportion of $\mcG \in \mbG'_{n,R}$ which have property $P$
approaches $c$ as $n \to \infty$ and hence also the proportion of $\mcG \in \mbG_{n,R}$
which have property $P$ approaches $c$ as $n \to \infty$.
To summarize the argument, we have:\\
{\bf Conclusion.} {\em
To prove that the proportion of $\mcG \in \mbG_{n,R}$ which have a property $P$ 
approaches $c$ as $n \to \infty$, independently of the number of vertices with degree $R-2$ in $\mcG$,
it suffices to prove that for every degree sequence $\bar{d}$
which satisfies~(\ref{conditions on degree sequence}) and such that 
$\sum_{i=1}^n d_i$ is even, the proportion of $\mcG \in \mbG_{n,R}(\bar{d})$
which has $P$, $|\mbP_n(\bar{d})|\big/|\mbG_{n,R}(\bar{d})|$, 
approaches $c$ as $n \to \infty$ and 
$\big| |\mbP_n(\bar{d})|\big/|\mbG_{n,R}(\bar{d})| \ - \ c \big|$ is bounded
by a function which tends to zero and depends only on $P$ and $R$.
}

\subsection{Proof of Theorem~\ref{theorem about typical structure}}

\noindent
The lemmas of this section prove the different parts of Theorem~\ref{theorem about typical structure}.
Let $R \geq 2$ be an integer.
Until the proof of Lemma~\ref{distribution of pairs of vertices of degree R-1 on a given distance}
is finished, we assume that
$W_1, \ldots, W_n$ are disjoint sets, $d_i = |W_i|$,
$\sum_{i=1}^n d_i = 2m$ and that~(\ref{conditions on degree sequence}) holds.
This implies that 
$$2m \sim Rn.$$
As we observed in Sections~\ref{configurations and multigraphs}
and~\ref{configurations and simple graphs},
if $\bar{d} = (d_1, \ldots, d_n)$, then
$\mbG_{n,R}(\bar{d})$ can be identified with 
$$\{\mr{Graph}(\mcC) : \mcC \in \mbC'(W_1, \ldots, W_n)\}.$$
By Fact~\ref{proportion of configurations whose images are graphs}
and the discussion in Section~\ref{graphs with a given degree sequence}, 
it follows that
in order to prove that the proportion of $\mcG \in \mbG_{n,R}$ with property $P$
approaches $c$ as $n \to \infty$, it now suffices to prove that 
$$\Bigg| \ \frac{\big|\big\{\mcC \in \mbC'(W_1, \ldots, W_n) : \mr{Graph}(\mcC) \text{ has } P\big\}\big|}
{|\mbC'(W_1, \ldots, W_n)|} \ - \ c \ \Bigg| \ = \  o(1) \quad \text{ as } n \to \infty,$$
where the bound $o( \ )$ depends only on $P$ and $R$.

\begin{lem}\label{small cycles are far from vertices with small degree}
For all integers $p \geq 1$ and $q \geq 3$, the proportion of $\mcG \in \mbG_{n,R}$
which have a vertex $v$ of degree less than $R$ and a $p$-path from $v$ to a
vertex in a $q$-cycle approaches 0 as $n \to \infty$.
It follows that the proportion of $\mcG \in \mbG_{n,R}$ which have a vertex $v$ of 
degree less than $R$ such that $v$ belongs to a $q$-cycle approaches $0$ as $n \to \infty$.
\end{lem}

\noindent
{\bf Proof.}
Let $p \geq 1$ and $q \geq 3$.
If $\mcC \in \mbC'(W_1, \ldots, W_n)$ and $\mr{Graph}(\mcC)$ has 
a vertex $W_i$ with degree less than $R$ and a $p$-path from $W_i$ to a vertex in a $q$-cycle,
then $|W_i| < R$ and there are $W_j$, a $p$-path $e_1, \ldots, e_p$ (of $\mcC$) and 
a $q$-cycle $e'_1, \ldots, e'_q$ (of $\mcC$) such that $e_1 \cap W_i \neq \es$,
$e_p \cap W_j \neq \es$ and $e'_1 \cap W_j \neq \es$.
If $W_i$ belongs to a $q$-cycle of $\mr{Graph}(\mcC)$,
then let $W_j = W_i$ and note that a $p$-path as $e_1, \ldots, e_p$ of $\mcC$ as above need not exist,
but the argument below (with $p = 0$ and finding at most $q-1$ other 
$F$-classes than $W_i$) shows that the corresponding quotient (below)
still approaches 0 as $n \to \infty$; this will prove the second statement of the lemma.
Therefore it suffices to prove that the proportion of $\mcC \in \mbC(W_1, \ldots, W_n)$
such that for some $W_i$ with $|W_i| < R$ there are $W_j$, a $p$-path
$e_1, \ldots, e_p$ and 
a $q$-cycle $e'_1, \ldots, e'_q$ such that $e_1 \cap W_i \neq \es$,
$e_p \cap W_j \neq \es$ and $e'_1 \cap W_j \neq \es$
approaches 0 as $n \to \infty$.

By assumption~(\ref{conditions on degree sequence}) 
we can choose $W_i$ with cardinality less than $R$ in at most 
$n^{1/4} + \sqrt{(R + \delta(n))n}$ ways, where $\delta(n) \to 0$ as $n\to\infty$.
A $p$-path $e_1, \ldots, e_p$ and a $q$-cycle $e'_1, \ldots, e'_q$
such that $e_p \cap W_j \neq \es$ and $e'_1 \cap W_j \neq \es$ for some $j$
can intersect at most $p - 1 + q$ different $F$-classes other than $W_i$.
Hence the at most $p+q-1$ $F$-classes other than $W_i$ 
which are going to include the union of the $p$-path and
$q$-cycle as above can be chosen in at most
$n^{p+q-1}$ ways. 
Then the $2(p+q)$ elements which are going to form the union of the
edges $e_1, \ldots, e_p, e'_1, \ldots, e'_q$ can be chosen in no more than
$\big((p+q)R\big)^{2(p+q)}$ ways. Then a complete matching on these
$2(p+q)$ elements (that is, nonintersecting edges $e_1, \ldots, e_p, e'_1, \ldots, e'_q$)
can be chosen in $\mr{M}(2(p+q))$ ways.
Finally, a complete matching on the remaining $2m - (2p + 2q)$ elements can
be chosen in $\mr{M}(2m - (2p + 2q))$ ways.
Therefore, there is a constant $\alpha > 0$ depending only on $p$, $q$ and $R$
such that, for all sufficiently large $n$, 
the proportion of $\mcC \in \mbC(W_1, \ldots, W_n)$ for which
there are $W_i$ with $|W_i| < R$, $W_j$, a $p$-path $e_1, \ldots, e_p$ and 
a $q$-cycle $e'_1, \ldots, e'_q$ such that $e_1 \cap W_i \neq \es$,
$e_p \cap W_j \neq \es$ and $e'_1 \cap W_j \neq \es$
is at most
\begin{align*}
\frac{\alpha \ \sqrt{n} \ n^{p + q - 1} \ \mr{M}(2m - 2(p + q))}
{\mr{M}(2m)} \
&\sim \ 
\frac{\alpha \ n^{p + q - 1/2}}{(2m)^{p+q}} \ \sim \ 
\frac{\alpha \ n^{p + q - 1/2}}{(Rn)^{p+q}} \\ 
&= \
\frac{\alpha \ n^{-1/2}}{R^{p+q}} \ \to \ 0 \quad \text{ as } \ n \to \infty,
\end{align*}
where~(\ref{approximation of quotient of matchings}) was used in the first
asymptotic identity. 
\hfill $\square$
\\

\noindent
The next lemma is a well known result in the case of $R$-regular graphs \cite{Wor99}, and
can be proved in a similar way in that case.

\begin{lem}\label{no subgraph with more edges than vertices}
If $\mcH$ is a graph with more edges than vertices, then 
the proportion of $\mcG \in \mbG_{n,R}$ which have a subgraph that is
isomorphic to $\mcH$ approaches 0 as $n \to \infty$.
\end{lem}

\noindent
{\bf Proof.}
Let $\mcH$ be a graph with $p > 0$ vertices and $q > p$ edges.
If $\mcC \in \mbC(W_1, \ldots, W_n)$ (so in particular if $\mcC \in \mbC'(W_1, \ldots, W_n)$)
and $\mr{Graph}(\mcC)$ has a subgraph $\mcH'$ which is isomorphic to $\mcH$,
then there are $W_{i_1}, \ldots, W_{i_p}$ and $p+1$ edges of $\mcC$
(i.e. $E$-classes) $e_1, \ldots, e_{p+1}$ such that $e_j \subseteq W_{i_1} \cup \ldots \cup W_{i_p}$
for $j = 1, \ldots, p+1$.
Hence it is sufficient to prove that the proportion of $\mcC \in \mbC(W_1, \ldots, W_n)$
such that there are $W_{i_1}, \ldots, W_{i_p}$ and $p+1$ edges of $\mcC$,
$e_1, \ldots, e_{p+1}$, such that $e_j \subseteq W_{i_1} \cup \ldots \cup W_{i_p}$
for $j = 1, \ldots, p+1$ approaches 0 as $n \to \infty$.

We can choose $W_{i_1}, \ldots, W_{i_p}$ in at most $n^p$ ways,
and then choose $p+1$ disjoint 2-subsets of $\bigcup_{j=1}^p W_{i_j}$ in at most
$\alpha = \binom{Rp}{2}\binom{Rp-2}{2} \cdots \binom{Rp-2p}{2}$ ways.
The remaining $2m - 2(p+1)$ elements can be completely matched in 
$\mr{M}(2m - 2(p + 1))$ ways.
Consequently, the proportion of $\mcC \in \mbC(W_1, \ldots, W_n)$ 
such that there are $W_{i_1}, \ldots, W_{i_p}$ and $p+1$ edges of $\mcC$,
$e_1, \ldots, e_{p+1}$, such that $e_j \subseteq W_{i_1} \cup \ldots \cup W_{i_p}$
for $j = 1, \ldots, p+1$ is at most 
\begin{align*}
\frac{\alpha \ n^p \ \mr{M}(2m - 2(p + 1))}{\mr{M}(2m)} \ \sim \
\frac{\alpha \ n^p}{(2m)^{p + 1}} \ \sim \ 
\frac{\alpha n^p}{(Rn)^{p+1}} \ = \ 
\frac{\alpha}{R^{p+1} n} \ \to \ 0 \ \text{ as } \ n \to \infty. \ \square
\end{align*}

\begin{rem}\label{remark about small cycles near each other}{\rm
Observe that Lemma~\ref{no subgraph with more edges than vertices} implies the following:
For all integers $p \geq 0$ and $q \geq 3$,
the proportion of $\mcG \in \mbG_{n,R}$ for which there are $3 \leq i,j \leq q$,
an $i$-cycle and a different $j$-cycle within distance $p$ of each other,
approaches 0 as $n \to \infty$.
}\end{rem}

\begin{lem}\label{3 vertices of degree less than R close to each other}
For every integer $p > 0$, the proportion of $\mcG \in \mbG_{n,R}$ that have three distinct vertices
$v_1, v_2, v_3$ with degree less than $R$ such that for all distinct $i,j \in \{1,2,3\}$ there is a path of 
length at most $p$ from $v_i$ to $v_j$, approaches 0 as $n \to \infty$.
\end{lem}

\noindent
{\bf Proof.}
Let $p$ be a positive integer.
It suffices to prove that, for every choice of positive integers 
$p_1, p_2 \leq p$, the proportion of $\mcC \in \mbC(W_1, \ldots, W_n)$ with the
following property approaches 0 as $n \to \infty$:
\begin{itemize}
\item[$(*)$] There are $W_{i_1}, W_{i_2}, W_{i_3}$ with 
$|W_{i_1}|, |W_{i_2}|, |W_{i_3}| < R$ 
such that there is a $p_1$-path with endpoints in
$W_{i_3}$ and $W_{i_1}$ and a $p_2$-path with endpoints in $W_{i_3}$ and $W_{i_2}$.
\end{itemize}
Let $c > R$. By Assumption~(\ref{conditions on degree sequence}), for all large enough $n$
there are at most $\sqrt{cn}$ $F$-classes $W_i$ with cardinality less than $R$,
so $W_{i_1}, W_{i_2}, W_{i_3}$ can be choosen in at most $(cn)^{3/2}$ ways.
The other $p_1 + p_2 - 2$ $F$-classes which the $p_1$-path and $p_2$-path are going to intersect
can be chosen in at most $n^{p_1 + p_2 - 2}$ ways. 
The union of all edges in the two paths will contain $2(p_1 + p_2)$ elements,
which can be chosen in no more than $\big((p_1 + p_2 + 1)R\big)^{2(p_1 + p_2)}$ ways, 
since they belong to the already
chosen $F$-classes. A complete matching (forming two paths) on these elements can be
chosen in no more than $\mr{M}(2(p_1 + p_2))$ ways, and a complete matching on the remaining
$2m - 2(p_1 + p_2)$ elements can be chosen in $\mr{M}(2m - 2(p_1 + p_2))$ ways.
Hence, for some constant $\alpha$ depending only on $p_1$, $p_2$ and $R$,
the proportion of $\mcC \in \mbC(W_1, \ldots, W_n)$ that satisfy $(*)$ is at most
\begin{align*}
\frac{\alpha n^{3/2} n^{p_1 + p_2 - 2} \mr{M}(2m - 2(p_1 + p_2))}{\mr{M}(2m)} \ 
&\sim \ 
\frac{\alpha n^{p_1 + p_2 - 1/2}}{(2m)^{p_1 + p_2}} \ \sim \ 
\frac{\alpha n^{p_1 + p_2 - 1/2}}{(Rn)^{p_1 + p_2}} \\ 
&= \ 
\frac{\alpha n^{-1/2}}{R^{p_1 + p_2}} \ \to \ 0 \quad \text{ as } \ n \to \infty. \quad \qquad \square
\end{align*}

\begin{lem}\label{elements of degree less than R-1 are far away}
For every integer $p > 0$, the proportion of $\mcG \in \mbG_{n,R}$
which have vertices $v$ and $w$ such that $\deg_{\mcG}(v) \leq R-2$, $\deg_{\mcG}(w) \leq R-1$ and 
a $p$-path from $v$ to $w$ approaches 0 as $n \to \infty$.
\end{lem}

\noindent
{\bf Proof.}
Let $p > 0$.
It suffices to prove that the proportion of $\mcC \in \mbC(W_1, \ldots, W_n)$
which have $F$-classes $W_{i_1}$, $W_{i_2}$ such that $|W_{i_1}| \leq R-2$ and $|W_{i_2}| \leq R-1$
and a $p$-path starting in $W_{i_1}$ and ending in $W_{i_2}$ approaches 0 as $n \to \infty$.
By Assumption~(\ref{conditions on degree sequence}), 
one can choose $W_{i_1}$ and $W_{i_2}$ such that
$|W_{i_1}| \leq R-2$ and $|W_{i_2}| \leq R-1$ in at most 
$n^{1/4}\cdot (\sqrt{(R + \delta(n))n} + n^{1/4}) = \sqrt{R + \delta(n)} \cdot n^{3/4} + \sqrt{n}$ ways.
The $p-1$ $F$-classes other than $W_{i_1}$ and $W_{i_2}$ 
which the $p$-path is going to intersect can be chosen in
at most $n^{p-1}$ ways. The number of ways in which $2p$ elements from $W_{i_1}$, $W_{i_2}$ and
the other chosen $F$-classes can be chosen is bounded by a constant depending only
on $p$ and $R$; the same is true for the number of ways of chosing a $p$-path from these elements.
Finally, a complete matching on the remaining $2m - 2p$ elements can be chosen in $\mr{M}(2m - 2p)$ ways.
So the proportion of $\mcC \in \mbC(W_1, \ldots, W_n)$ with the described property is,
for some constant $\alpha$ depending only on $p$ and $R$, at most
$$\frac{\alpha n^{3/4} \ n^{p-1} \ \mr{M}(2m - 2p)}{\mr{M}(2m)} \ \sim \ 
\frac{\alpha n^{p - 1/4}}{(2m)^p} \ \sim \ \frac{\alpha n^{p - 1/4}}{(Rn)^p} \ = \ 
\frac{\alpha n^{-1/4}}{R^p} \ \to \ 0 \quad \text{ as } \ n \to \infty. \ \hfill \square$$

\begin{lem}\label{distribution of pairs of vertices of degree R-1 on a given distance}
Let $k \geq 3$ be an integer and let $r_3, \ldots, r_k$, $s_1, \ldots, s_k \in \mbbN$.
Then the proportion of $\mcG \in \mbG_{n,R}$ which, for $p = 3, \ldots, k$, have exactly 
$r_p$ $p$-cycles and, for $p = 1, \ldots, k$, have exactly $s_p$ $p$-paths with both endpoints
of degree $R-1$, approaches
$$\Bigg(\prod_{p = 3}^k \frac{(\lambda_p)^{r_p} \ e^{-\lambda_p}}{r_p!}\Bigg)
\Bigg(\prod_{p = 1}^k \frac{(\mu_p)^{s_p} \ e^{-\mu_p}}{s_p!}\Bigg)
\quad \text{ as } \ n \to \infty,$$
where $\lambda_p = \frac{(R-1)^p}{2p}$ and $\mu_p = \frac{(R-1)^{p+1}}{2}$.
In other words, if $X_p$ is the number of $p$-cycles and $Y_p$ is the number of 
$p$-paths with both endpoints of degree $R-1$, then the random variables $X_p$, $p = 3, \ldots, k$
and $Y_p$, $p = 1, \ldots, k$, are asymptotically independent Poisson variables with
mean $\lambda_p$ and $\mu_p$, respectively.
Moreover, the random variables $X_p$, $p = 3, \ldots, k$, and $Y_p$, $p = 1, \ldots, k$,
are asymptotically independent of the (asymptotically Poisson distributed)
number of vertices with degree $R-2$.
\end{lem}

\noindent
{\bf Proof.} 
Let $k \geq 3$.
The last statement of the lemma will follow from the argument below,
because the only thing regarding the
number of vertices with degree $R-2$, or in the context of configurations, 
the number of $F$-classes of cardinality $R-2$, is that this number is
at most $n^{1/4}$; so the limit to be proved is independent of the
number of vertices ($F$-classes) with degree (cardinality) $R-2$. 
For $p = 1, \ldots, k$, let the random variable $X_p$ be the number of 
$p$-cycles of a configuration $\mcC \in \mbC(W_1, \ldots, W_n)$
and let the random variable $Y_p$ be the number of $p$-paths (of $\mcC$) with both endpoints 
in $F$-classes of cardinality $R-1$.
For $p = 1, \ldots, k$, let 
$$\lambda_p = \frac{(R-1)^p}{2p} \quad \text{ and } \quad 
\mu_p = \frac{(R-1)^{p+1}}{2}.$$
We will prove:
\\

\noindent
{\em Claim.} $X_1, \ldots, X_k, Y_1, \ldots, Y_k$ are asymptotically independent
Poisson variables with mean $\lambda_1, \ldots, \lambda_k, \mu_1, \ldots, \mu_k$, respectively.
\\

\noindent
By Fact~\ref{proportion of configurations whose images are graphs}~(i), 
the probability that $X_1 = X_2 = 0$ approaches 
$\exp(-\lambda_1 - \lambda_2)$ as $n \to \infty$, and therefore 
(using part (ii) of the same fact) the lemma follows
from the claim. Actually the argument that follows proves 
Fact~\ref{proportion of configurations whose images are graphs}~(i),
but also considers the random variables $Y_1, \ldots, Y_k$ of the lemma, 
which we must take into account when proving a logical limit law.

Recall the notation $(x)_i = x(x-1) \ldots (x-i+1)$.
Note that the random variable
$$Z = (X_1)_{r_1}(X_2)_{r_2} \ldots (X_k)_{r_k} (Y_1)_{s_1}(Y_2)_{s_2} \ldots (Y_k)_{s_k}$$
is the number of ordered $2k$-tuples without repetition consisting of (from left to right) 
$r_1$ 1-cycles, $r_2$ 2-cycles, $\ldots$, $r_k$ $k$-cycles, 
$s_1$ 1-paths with both endpoints in $F$-classes of cardinality $R-1$, $\ldots$ and
$s_k$ $k$-paths with both endpoints in $F$-classes of cardinality $R-1$.
By Theorem~1.23 in~\cite{Bol01}, to prove the claim it suffices to prove that
\begin{equation}\label{convergence of factorial moments}
\Bigg| \mbbE(Z) \ - \ 
\prod_{p = 1}^k (\lambda_p)^{r_p} (\mu_p)^{s_p} \Bigg| \ = \ o(1) \quad \text{ as } n \to \infty,
\end{equation}
where $\mbbE( \ )$ denotes the expected value and the bound `$o( \ )$' depends only
on $k$, $R$, $r_1, \ldots, r_k$ and $s_1, \ldots, s_k$.
Indeed, in the argument we never use any other properties of the degree sequence
$\bar{d} = (d_1, \ldots, d_n)$ than those stated in~({\ref{conditions on degree sequence}),
which hold for all degree sequences that we consider.
Let $Z'$ be the number of ordered $2k$-tuples $(A_1, \ldots, A_{2k})$
without repetition as counted by $Z$ above,
but with the extra condition that there is {\em no} $W_i$ such that for two distinct entries 
$A_j$, $A_{j'}$ ($j' \neq j$) and edges $e_j \in A_j$ and $e_{j'} \in A_{j'}$, 
$e_j \cap W_i \neq \es$ and $e_{j'} \cap W_i \neq \es$.
Let $Z'' = Z - Z'$, so $Z = Z' + Z''$ and $\mbbE(Z) = \mbbE(Z') + \mbbE(Z'')$.
Let $C_p(n)$ be the number of ways to choose a $p$-cycle (a matching on a $2p$-subset of 
$W_1 \cup \ldots \cup W_n$ which forms a $p$-cycle), so we have
\begin{equation}\label{asymptotic for C-p}
C_p(n) \ \leq \ \frac{\big(R(R-1)n\big)^p}{2p}.
\end{equation}
Then let $P_p(n)$ be the number of ways to choose a $p$-path with both
endpoints in $F$-classes with cardinality $R-1$, so by 
assumption~(\ref{conditions on degree sequence}),
\begin{align}\label{asymptotic for P-p}
P_p(n) \ 
&\leq \ \frac{(R-1)^2(R + \delta_n)n \ \big(R(R-1)n\big)^{p-1}}{2} \\ 
&\sim \ \frac{(R-1)\big(R(R-1)n\big)^p}{2} \quad \text{ as } \ n \to \infty. \nonumber
\end{align}
For every $A \subset \bigcup_{i=1}^n W_i$, let 
$C_p(n, A)$ and $P_p(n, A)$ be defined 
as $C_p(n)$ and $P_p(n)$ but with the extra condition that
no choosen edge ($E$-class) has non-empty intersection with $A$.
Let $t = \sum_{p=1}^k (pr_p + ps_p)$ and let
\begin{align*}
C_p(n,t) \ &= \ \min\big\{C_p(n, A) : A \subseteq \bigcup_{i=1}^n W_i, \ |A| = 2t\big\}, \ \text{ and} \\
P_p(n,t) \ &= \ \min\big\{C_p(n, A) : A \subseteq \bigcup_{i=1}^n W_i, \ |A| = 2t\big\}.
\end{align*}
Let $f(n) \ = \ n \ - \ \sqrt{(R + \delta_n)n} \ - \ n^{1/4} \ - \ 2t$.
From~(\ref{conditions on degree sequence}) we get
\begin{align}
C_p(n,t) \ &\geq \ \frac{\big( R(R-1)f(n) \big)^p}{2p} \ \sim \ 
\frac{\big( R(R-1)n \big)^p}{2p} \quad \text{ and} 
\label{asymptotic for C-p(n,t)} \\
P_p(n,t) \ &\geq \ \frac{(R-1)^2 \big(\sqrt{(R - \delta_n)n} \ - \ 2t \big)^2}{2} \ \big( R(R-1)f(n)\big)^{p-1} 
\label{asymptotic for P-p(n,t)} \\ 
&\sim \ 
\frac{(R-1)\big(R(R-1)n\big)^p}{2}  \quad \text{ as } \ n \to \infty. \nonumber
\end{align}
We have
$$\prod_{p=1}^k C_p(n,t)^{r_p} P_p(n,t)^{s_p}  \ \leq \ Z' \ \leq \ 
\prod_{p=1}^k C_p(n)^{r_p} P_p(n)^{s_p},$$
so by~(\ref{asymptotic for C-p}),~(\ref{asymptotic for P-p}),~(\ref{asymptotic for C-p(n,t)}) 
and~(\ref{asymptotic for P-p(n,t)}) we get
\begin{equation}\label{asymptotic for Z'}
Z' \ \sim \ \prod_{p=1}^k \Bigg( \frac{\big(R(R-1)n\big)^p}{2p} \Bigg)^{r_p} 
\Bigg( \frac{(R-1)\big(R(R-1)n\big)^p}{2} \Bigg)^{s_p}.
\end{equation}
For every tuple $(A_1, \ldots, A_{2k})$ that is being counted by $Z'$, we
have $\Big|\bigcup_{p=1}^{2k} A_p\Big| = 2t$ and for every 
$A \subseteq \bigcup_{i=1}^n W_i$ with $|A| = 2t$ and every complete matching on $A$, 
the proportion of configurations in $\mbC(W_1, \ldots, W_n)$ which have this matching
on A is
$$\frac{\mr{M}(2m - 2t)}{\mr{M}(2m)} \ \sim \ (2m)^{-t} \ \sim \ (Rn)^{-t}.$$
From this and~(\ref{asymptotic for Z'}) and recalling that $t = \sum_{p=1}^k (pr_p + ps_p)$, we get
\begin{align*}
\mbbE(Z') \ 
&\sim \ (Rn)^{-\sum_{p=1}^k (pr_p + ps_p)} \ 
\prod_{p=1}^k \Bigg( \frac{\big(R(R-1)n\big)^p}{2p} \Bigg)^{r_p} 
\Bigg( \frac{(R-1)\big(R(R-1)n\big)^p}{2} \Bigg)^{s_p} \\
&= \ \prod_{p=1}^k \Bigg( \frac{(R-1)^p}{2p} \Bigg)^{r_p} 
\Bigg( \frac{(R-1)^{p+1}}{2} \Bigg)^{s_p} \ = \ 
\prod_{p=1}^k (\lambda_p)^{r_p} (\mu_p)^{s_p}.
\end{align*}
By~(\ref{asymptotic for C-p}),~(\ref{asymptotic for P-p}),~(\ref{asymptotic for C-p(n,t)}) 
and~(\ref{asymptotic for P-p(n,t)}), 
$\big|\mbbE(Z') \ - \ \prod_{p=1}^k (\lambda_p)^{r_p} (\mu_p)^{s_p}\big| =~o(1)$
as $n \to \infty$ where the bound depends only on
$R$, $k$, $r_1, \ldots, r_k$ and $s_1, \ldots, s_k$.

Since $Z = Z' + Z''$ it now suffices to prove that $\lim_{n\to\infty}\mbbE(Z'') =~o(1)$,
where the bound depends only on $R$, $k$, $r_1, \ldots, r_k$ and $s_1, \ldots, s_k$,
because this together with the conclusion above implies 
the statement of~(\ref{convergence of factorial moments}).
But $Z'' > 0$ means that at least one of the following conditions holds:
\begin{itemize}
\item[(a)] For some $q \leq t$ there are $W_{i_1}, \ldots, W_{i_q}$ and edges $e_1, \ldots, e_{q+1}$ such that
$e_j \subseteq W_ {i_1} \cup \ldots \cup W_{i_q}$ for $j = 1, \ldots, q+1$.

\item[(b)] For some $q, q' \leq t$ and some $W_i$ with $|W_i| < R$ there are $W_j$, a $q$-path $e_1, \ldots, e_q$
and $q'$-cycle $e'_1, \ldots, e'_{q'}$ such that $e_1 \cap W_i \neq \es$,
$e_q \cap W_j \neq \es$ and $e'_1 \cap W_j \neq \es$.

\item[(c)] There are distinct $W_{i_1}, W_{i_2}, W_{i_3}$
with $|W_{i_1}|, |W_{i_2}|, |W_{i_3}| < R$ and for every pair $j, j' \in \{1,2,3\}$ such that $j \neq j'$
a path of length at most $t$ from some $a \in W_{i_{j}}$ to some $b \in W_{i_{j'}}$.
\end{itemize}
From the proofs of Lemmas~\ref{small cycles are far from vertices with small degree},
~\ref{no subgraph with more edges than vertices} and~\ref{3 vertices of degree less than R close to each other},
it follows that the proportion of $\mcC \in \mbC(W_1, \ldots, W_n)$ satisfying any of~(a),~(b) or~(c) 
approaches 0 as $n \to \infty$, and in each case the convergence is bounded by a function depending
only on $R$, $k$, $r_1, \ldots, r_k$ and $s_1, \ldots, s_k$. 
Hence $\mbbE(Z'') = o(1)$ as $n \to \infty$, where the bound depends only on
$R$, $k$, $r_1, \ldots, r_k$ and $s_1, \ldots, s_k$.
\hfill $\square$

\begin{lem}\label{no small components}
Suppose that $R \geq 3$. For every integer $p$, the proportion of $\mcG \in \mbG_{n,R}$ such that
every connected component of $\mcG$ has at least $p$ vertices approaches 1 as $n \to \infty$.
\end{lem}

\noindent
{\bf Proof.}
Let $p$ be any positive integer.
By Lemmas~\ref{no subgraph with more edges than vertices},
~\ref{3 vertices of degree less than R close to each other},
~\ref{elements of degree less than R-1 are far away}
and Theorem~\ref{theorem about typical degree distribution},
almost all $\mcG \in \mbG_R$ have the following properties:
\begin{itemize}
\item[\textbullet] \ There is no subgraph with $p$ vertices and more than $p$ edges.
\item[\textbullet] \  There do not exist distinct vertices $v_i$, $i = 1,2,3$, all three with
degree $R-1$ such that for all distinct $i,j \in \{1,2,3\}$ there is a path of length
at most $p$ from $v_i$ to $v_j$.
\item[\textbullet] \ There do not exist distinct vertices $v$ and $w$ both of degree
$R-2$ and a path of length at most $p$ from $v$ to $w$.
\item[\textbullet] \ No vertex has degree less than $R-2$.
\end{itemize}
Therefore it suffices to show that if $\mcG \in \mbG_R$ has a connected component
with exactly $p$ vertices, then one of the above properties fail for $\mcG$.
Recall the assumption that $R \geq 3$ and suppose that $\mcH \subseteq \mcG$ 
is a connected component with exactly $p$ vertices.
If $p = 1$ then the unique vertex in $\mcH$ has degree $0 < R-2$, so $\mcG$ does not have
the last property above.
If $p = 2$ then the two vertices of $\mcH$ have degree $1 \leq R-2$ in $\mcG$, so 
the third or fourth property fails.
If $p = 3$ then the three verties of $\mcH$ have degree at most $2 \leq R-1$ in $\mcG$, so
$\mcG$ does not have the second property.
Now suppoes that $p \geq 4$.
If the last three properties hold, then the number of edges in $\mcH$ is at least
$$\frac{(R-2) \ + \ (R-1) \ + \ (p-2)R}{2} \ = \ \frac{pR - 3}{2} \ \geq \ p\frac{3}{2} - \frac{3}{2}
\ > \ p,$$
where the last inequality holds for all $p \geq 4$ (by induction or differentiation).
Hence $\mcH$ has more than $p$ edges, so $\mcG$ does not have the first property.
\hfill $\square$
\\

\noindent
Note that 
Lemma~\ref{no small components} is false for $R=2$, because if $R=2$ then every $3$-cycle (say)
is a connected component and the proportion of $\mcG \in \mbG_{n,2}$ with at least
one $3$-cycle converges to a positive number, by 
Lemma~\ref{distribution of pairs of vertices of degree R-1 on a given distance}, as $n \to \infty$.
In the case $R \geq 5$ Lemma~\ref{no small components}
is a consequence of Corollary~\ref{corollary to main theorem}.

\section{First-order limit laws}\label{First-order limit laws}

\noindent
Recall that $\mbG_{n,R}$ denotes the set of undirected graphs with 
vertices $1, \ldots, n$ such that every vertex has degree at most $R$.
For every $n$ and $\mcG \in \mbG_{n,R}$ we let $[\mcG]$ denote the equivalence class to which it belongs
with respect to the isomorphism relation on $\mbG_{n,R}$.
Note that for every first-order sentence $\varphi$ and every $\mcG \in \mbG_{n,R}$,
we have $\mcG \models \varphi$ if and only if $\mcH \models \varphi$ for every $\mcH \in [\mcG]$.

\begin{theor}\label{the limit law}
Let $R$ be a non-negative integer. \\
(i) Suppose that $R \geq 0$.
For every first order sentence $\varphi$ in the language of graphs, there is $c \in [0,1]$
such that
$$\lim_{n\to\infty} \ 
\frac{\big|\{\mcG \in \mbG_{n,R} : \mcG \models \varphi \}\big|}{\big|\mbG_{n,R}\big|} \ = \ c.$$
(ii) Suppose that $0 \leq R \leq 1$ or $R \geq 5$.
For every first-order sentence $\varphi$ in the language of graphs, there is $c \in [0,1]$
such that
$$\lim_{n\to\infty} \ 
\frac{\big|\{[\mcG] : \mcG \in \mbG_{n,R} \text{ and } \mcG \models \varphi\}\big|}
{\big|\{[\mcG] : \mcG \in \mbG_{n,R}\}\big|} \ = \ 
\lim_{n\to\infty} \ \frac{\big|\{\mcG \in \mbG_{n,R} : \mcG \models \varphi \}\big|}{\big|\mbG_{n,R}\big|} 
\ = \ c.$$
\end{theor}

\noindent
In other words, for every $R \geq 0$, finite graphs with maximum degree $R$ satisfy a
labelled limit law for first-order logic. If $0 \leq R \leq 1$ or $R \geq 5$ then we also
have an unlabelled limit law.
If $R=0$ or $R=1$ then we have a zero-one law in both the labelled an unlabelled case
(as the proof below will show), that is, the number $c$ in the theorem is
either 0 or 1 for every $\varphi$.
If $R \geq 2$ then we do not have a zero-one law, because (for example) the non-existence
of a vertex with degree $R-2$, which can be expressed with a sentence in first-order logic, 
holds with asymptotic probability $e^{-(R-1)}$, by 
Theorem~\ref{theorem about typical degree distribution}.
This paper does not settle the question of whether or not $\mbG_{n,R}$ has an unlabelled limit
law for $R = 2, 3, 4$.

\subsection{Proof of Theorem~\ref{the limit law}}\label{proof of the limit law}

If $R = 0$ then $|\mbG_{n,R}| = 1$ for every $n$, so 
the theorem is trivial in this case.
If $R = 1$ then there are only two isomorphism types of connected components, singletons
and two vertices connected to each other, so by Example~7.15 in \cite{Com87}
the theorem follows, both in the labelled and unlabelled case.
(In this simple case one can of course also argue directly, without reference to \cite{Com87}
which has much wider applicability.)

For the rest of the proof we assume that $R \geq 2$.
We will prove that for every first-order sentence $\varphi$, the proportion
$|\{\mcG \in \mbG_{n,R} : \mcG \models \varphi\}| \big/ |\mbG_{n,R}|$ converges as $n \to \infty$;
so we get part~(i) of the theorem for all $R \geq 2$.
By Corollary~\ref{corollary to main theorem} and 
Theorem~\ref{rigidity implies that labelled and unlabelled probabilities coincide} 
we then get part~(ii) for $R \geq 5$, because
in Theorem~\ref{rigidity implies that labelled and unlabelled probabilities coincide} 
we can, for any sentence $\varphi$, let $\mbH_n$ be the set of 
$\mcG \in \mbG_{n,R}$ such that $\mcG \models \varphi$.

Every first-order sentence has a quantifier rank (also called quantifier depth) 
(see \cite{EF, Lib, Spe} for example) 
which is a non-negative integer.
Therefore it suffices to show that for every $k > 0$ and every first-order sentence $\varphi$ with
quantifier rank at most $k$, the quotient 
$|\{\mcG \in \mbG_{n,R} : \mcG \models \varphi \}| \big/ |\mbG_{n,R}|$ converges as $n \to \infty$.
So we fix an arbitrary integer $k \geq 3$.

\begin{defin}{\rm
Let $m = 5^{k+1}$. 
For all $q, r_3, \ldots, r_m, s_1, \ldots, s_m \in \{0, 1, \ldots, k\}$, define 
$$\mbX_n(q,r_3, \ldots, r_m, s_1, \ldots, s_m)$$
to be the set of all $\mcG \in \mbG_{n,R}$ such that the following holds:
\begin{enumerate}
\item There are no vertices with degree less that $R-2$,
\item If $q < k$ then there are exactly $q$ vertices with degree $R-2$ and if
$q = k$ then there are at least $k$ vertices with degree $R-2$.
\item There are at least $m$ vertices with degree $i$ for $i = R-1, R$. 
\item For $p = 3, \ldots, m$, if $r_p < k$ then there are exactly $r_p$ $p$-cycles and
if $r_p = k$ then there are at least $k$ $p$-cycles.
\item For $p = 1, \ldots, m$, if $s_p < k$ then there are exactly $s_p$ $p$-paths with
both endpoints of degree $R-1$ and if $s_p = k$ then there are at least $k$ $p$-paths with
both endpoints of degree $R-1$.
\item The distance is at least $5^{k+2}$ between 
\begin{itemize}
\item[(a)] any vertex with degree $R-2$ and any (other) vertex of degree at most $R-1$,
\item[(b)] any vertex with degree at most $R-1$ and any cycle of length at most $m$, 
\item[(c)] any two different cycles of length at most $m$, and
\item[(d)] any two different paths of length at most $m$ such that both endpoints of both
paths have degree $R-1$.
\end{itemize}
\item If $R \geq 3$ then every connected component has at least $m$ vertices.
\end{enumerate}
}\end{defin}

\noindent
Let 
\[
P_k(x,\mu) \ = \ 
\begin{cases}
\frac{\mu^x e^{-\mu}}{x!} \quad \text{if } \ x < k, \\
1 \ - \ \sum_{i=0}^{k-1} \frac{\mu^i e^{-\mu}}{i!} \quad {if } \ x \geq k,
\end{cases}
\]
and let $\lambda_p = \frac{(R-1)^p}{2p}$ and $\mu_p = \frac{(R-1)^{p+1}}{2}$.
By Theorems~\ref{theorem about typical degree distribution}
and~\ref{theorem about typical structure}
and Lemma~\ref{distribution of pairs of vertices of degree R-1 on a given distance},
\begin{align*}
\lim_{n\to\infty} \ &\frac{|\mbX_n(q,r_3, \ldots, r_m, s_1, \ldots, s_m)|}{|\mbG_{n,R}|} \\
= \ &P_k(q, R-1) \ 
\Bigg(\prod_{p = 3}^m P_k(r_p, \lambda_p)\Bigg)
\Bigg(\prod_{p = 1}^m P_k(s_p, \mu_p)\Bigg).
\end{align*}
Hence it is enough to prove that for all
$q, r_3, \ldots, r_m, s_1, \ldots, s_m \in \{0, 1, \ldots, k\}$,
if $\mcG, \mcH \in \mbX_n(q, r_3, \ldots, r_m, s_1, \ldots, s_m)$, 
then $\mcG$ and $\mcH$ satisfies exactly the same sentences of
quantifier rank at most $k$.
To show that $\mcG$ and $\mcH$ satisfies exactly the same sentences
of quantifier rank at most $k$ it suffices to prove that 
Duplicator has a winning strategy for the Ehrenfeucht-Fra\"{i}ss\'{e} game
in $k$ steps on $\mcG$ and $\mcH$ (e.g. \cite{EF}, Theorem~2.2.8, or
similar results in~\cite{Lib, Spe}).

So we fix $q, r_3, \ldots, r_m, s_1, \ldots, s_m \in \{0, 1, \ldots, k\}$ and let
$$\mcG, \mcH \ \in \ \mbX_n(q, r_3, \ldots, r_m, s_1, \ldots, s_m).$$
For any vertex $v$ of $\mcG$ and $l \in \mbbN$ we let $B_{\mcG}(v,l)$
be the set of vertices $w$ such that $\dist_{\mcG}(v,w) \leq l$; and similarly for $\mcH$.
Let a {\em Poisson object} of $\mcG$ or $\mcH$ denote any one of
\begin{itemize}
\item[\textbullet] \ a vertex with degree $R-2$, or
\item[\textbullet] \ a $p$-cycle where $p \leq m$, or
\item[\textbullet] \ a $p$-path with both endpoints of degree $R-1$ where $p \leq m$.
\end{itemize}

\begin{observation}\label{observation}{\rm
It follows from the definition of 
$\mbX_n(q, r_3, \ldots, r_m, s_1, \ldots, s_m)$ that 
if $v \in \{1, \ldots, n\}$, $B = B_{\mcG}(v,l)$, or $B = B_{\mcH}(v,l)$, and $l \leq 5^k$, then
$B$ contains at most one Poisson object.
}\end{observation}

\noindent
To prove that Duplicator has a winning strategy for the Ehrenfeucht-Fra\"{i}ss\'{e} game
in $k$ steps on $\mcG$ and $\mcH$ it suffices to prove the following result.

\begin{lem}
Suppose that $i < k$, $v_1, \ldots, v_i, w_1, \ldots, w_i \in \{1, \ldots, n\}$,
$v_{i+1} \in \{1, \ldots, n\}$ (or $w_{i+1} \in \{1, \ldots, n\}$) and that
\begin{equation*}\label{the ith isomorphism}
f_i : \mcG\bigg[\bigcup_{j=1}^i B_{\mcG}\big(v_i, 5^{k-i}\big)\bigg] \ \to \ 
\mcH\bigg[\bigcup_{j=1}^i B_{\mcH}\big(w_i, 5^{k-i}\big)\bigg]
\end{equation*}
is an isomorphism such that $f_i(v_j) = w_j$ for $j = 1, \ldots, i$.
Then there is $w_{i+1} \in \{1 ,\ldots, n\}$ (or $v_{i+1} \in \{1, \ldots, n\}$) 
such that the same statement holds with `$i+1$' in place of `$i$'.
\end{lem}

\noindent
{\bf Proof.}
Suppose that $i < k$, $v_1, \ldots, v_i, w_1, \ldots, w_i \in \{1, \ldots, n\}$ and that
$$f_i : \mcG\bigg[\bigcup_{j=1}^i B_{\mcG}\big(v_i, 5^{k-i}\big)\bigg] \ \to \ 
\mcH\bigg[\bigcup_{j=1}^i B_{\mcH}\big(w_i, 5^{k-i}\big)\bigg]$$
is an isomorphism such that $f_i(v_j) = w_j$ for $j = 1, \ldots, i$.
By symmetry it is enough to consider the case when Spoiler chooses $v_{i+1}$ in $\mcG$.

If
$B_{\mcG}\big(v_{i+1},5^{k-i-1}\big) \subseteq 
\bigcup_{j=1}^i B_{\mcG}\big(v_j, 5^{k-i} - 1\big)$
then let $w_{i+1} = f_i(v_{i+1})$ and let $f_{i+1}$ be the restriction of $f_i$
to $\bigcup_{j=1}^{i+1} B_{\mcG}\big(v_j, 5^{k-i-1}\big)$.

Now suppose that $B_{\mcG}\big(v_{i+1},5^{k-i-1}\big)$ is not 
included in $\bigcup_{j=1}^i B_{\mcG}\big(v_j, 5^{k-i} - 1\big)$.
Then there is $u \in B_{\mcG}\big(v_{i+1}, 5^{k-i-1}\big)$ such that the
distance (in $\mcG$) from $u$ to $v_j$ is at least $5^{k-i}$ for
all $j = 1, \ldots, i$.
It follows that for every $a \in B_{\mcG}\big(v_{i+1}, 5^{k-i-1}\big)$ and
every $j = 1, \ldots, i$, the distance from $a$ to $v_j$ is at least
$$5^{k-i} \ - \ \dist_{\mcG}(a,u) \ \geq \ 5^{k-i} \ - \ 2\cdot 5^{k-i-1} \ = \ 3\cdot 5^{k-i-1}.$$
Consequently,  for every $a \in B_{\mcG}\big(v_{i+1}, 5^{k-i-1}\big)$
and every $b \in \bigcup_{j=1}^i B_{\mcG}\big(v_j, 5^{k-i-1}\big)$,
$$\dist_{\mcG}(a,b) \ \geq \ 3\cdot 5^{k-i-1} \ - \ 5^{k-i-1} \ = \ 2\cdot 5^{k-i-1} \ \geq \ 2,$$
because $i < k$.
Since $\mcG, \mcH \in \mbX_n(q, r_3, \ldots, r_m, s_1, \ldots, s_m)$ and
$f_i$ is an isomorphism it follows from Observation~\ref{observation}, for large enough $n$,
that there are $w_{i+1}$ and an isomorphism 
$$f'_{i+1} : \mcG \Big[B_{\mcG}\big(v_{i+1}, 5^{k-i-1}\big) \Big] \ \to \ 
\mcH \Big[B_{\mcH}\big(w_{i+1}, 5^{k-i-1}\big) \Big]$$
such that $f'_{i+1}(v_{i+1}) = w_{i+1}$ and $\dist_{\mcH}(w_{i+1},u) \geq 2$ 
for every $u \in \bigcup_{j=1}^i B_{\mcH}\big(w_i, 5^{k-i}\big)$.
Then the desired isomorphism $f_{i+1}$ is obtained by letting
$f_{i+1}(v) = f_i(v)$ if \\
$v \in \bigcup_{j=1}^i B_{\mcG}\big(v_i, 5^{k-i}\big)$
and $f_{i+1}(v) = f'_{i+1}(v)$ if $v \in B_{\mcG}\big(v_{i+1}, 5^{k-i-1}\big)$.
\hfill $\square$

\bigskip

\noindent
{\bf Acknowledgements.}
I thank Svante Janson and Nicholas Wormald for helpful 
comments and discussions. In particular, Janson suggested the proof of
Theorem~\ref{theorem about typical degree distribution} which is presented here.
(My proof first obtained some rough estimates by counting arguments
and then applied a known asymptotic formula for graphs with given degree sequence,
which, however, involved rather long and tedious calculations.)
Also, thanks to the anonymous referee for careful reading and 
helping me to remove various mistakes.

\end{document}